\documentclass[12pt]{article}
\usepackage{amssymb}

\newcommand{\be}{\begin{equation}}
\newcommand{\ee}{\end{equation}}

\newcommand{\bea}{\begin{eqnarray}}
\newcommand{\eea}{\end{eqnarray}}

\newcommand{\eq}{\! & = & \!}
\newcommand{\eqq}{\!\! & = & \!\!}

\newcommand{\st}[2]{\rule[#1]{0mm}{#2}} 

\newcommand{\hst}[1]{\rule{#1}{0mm}}

\newcommand{\ty}{\hspace{0.04em}}
\newcommand{\tty}{\hspace{0.06em}}

\newcommand{\nn}{\nonumber}

\newcommand{\sect}[1]{\setcounter{equation}{0}\section{#1}}

\newcommand{\reff}[1]{(\ref{#1})}

\newcommand{\dis}{\displaystyle}

\newcommand{\Ga}{\Gamma}
\newcommand{\De}{\Delta}

\newcommand{\al}{\alpha}
\newcommand{\bet}{\beta}
\newcommand{\ga}{\gamma}
\newcommand{\de}{\delta}

\newcommand{\ve}{\varepsilon}

\newcommand{\lam}{\lambda}

\newcommand{\si}{\sigma}
\newcommand{\vp}{\varphi}

\newcommand{\la}{\langle}
\newcommand{\ra}{\rangle}
\newcommand{\rar}{\rightarrow}
\newcommand{\lra}{\longrightarrow}
\newcommand{\ti}{\times}
\newcommand{\op}{\oplus}

\newcommand{\ot}{\otimes}

\newcommand{\KK}{\mathbb K}

\newcommand{\ZZ}{\mathbb Z}

\newcommand{\CC}{\mathbb C}
\newcommand{\II}{\mathbb I}

\newcommand{\ohal}{{\textstyle\frac{1}{2}}}

\newcommand{\cd}{\cdot}

\newcommand{\comp}{\!\stackrel{\textstyle{\hst{0ex} \atop \circ}}{\hst{0ex}}\!}

\newcommand{\mc}{\multicolumn}
 
\newcommand{\ub}{\underbrace}
\newcommand{\ol}{\overline}
\newcommand{\spo}{{\mathfrak{spo}}(2n \ty|\ty 2m)}
\newcommand{\spott}{{\mathfrak{spo}}(2 \ty|\ty 2)}
\newcommand{\spob}{{\mathfrak{spo}}(b)}
\newcommand{\osp}{{\mathfrak{osp}}(2m \ty|\ty 2n)}
\newcommand{\gl}{{\mathfrak{gl}}}
\newcommand{\sli}{{\mathfrak{sl}}}
\newcommand{\fh}{{\mathfrak{h}}}
\newcommand{\Uq}{U_q({\mathfrak{spo}}(2n \ty|\ty 2m))}
\newcommand{\id}{\mbox{id}}
\newcommand{\ad}{\mbox{ad}}
\newcommand{\add}{\mbox{ad}\tty}
\newcommand{\adb}{\overline{\mbox{ad}}}
\newcommand{\addb}{\overline{\mbox{ad}}\tty}
\newcommand{\Rh}{\hat{R}}
\newcommand{\Cqi}{((C^q)^{-1})}
\newcommand{\Lgr}{\mbox{Lgr}}
\newcommand{\otb}{\:\overline{\otimes}\:}
\newcommand{\SP}{\mbox{SPO}_q(2n \ty|\ty 2m)}

\newcommand{\bt}{\tilde{b}}

\newcommand{\Vf}{\overline{V}_{\!\!4}}
\newcommand{\Ch}{{\mathbb C}[[h]]}
\newcommand{\Vh}{V[[h]]}
\newcommand{\otC}{\!\begin{array}[t]{c} \ot \\[-1.2ex]
                                    {\scriptstyle\CC} \end{array}\!} 
\newcommand{\otCh}{\!\begin{array}[t]{c} \ot \\[-1.2ex]
                                    {\scriptstyle\CC[[h]]} \end{array}\!} 
\newcommand{\su}{{\displaystyle{{}^s}}}
\newcommand{\au}{{\displaystyle{{}^a}}}
\newcommand{\vso}{\vspace{-1.0ex}}

\begin{document}

\begin{titlepage}

\newlength{\ppn}
\newcommand{\defboxn}[1]{\settowidth{\ppn}{#1}}
\defboxn{BONN--TH--2000--03}

\hspace*{\fill} \parbox{\ppn}{BONN--TH--2000--03 \\
                              April 2000        }

\vspace{18mm}

\begin{center}
{\LARGE\bf The $R$--matrix \\[0.4ex]
           of the symplecto--orthogonal \\[0.5ex]
           quantum superalgebra $\Uq$ \\[0.7ex]
           in the vector representation} \\
\vspace{15mm}
{\large M. Scheunert} \\
\vspace{1mm}
Physikalisches Institut der Universit\"{a}t Bonn \\
Nu{\ss}allee 12, D--53115 Bonn, Germany \\
\end{center}

\vspace{15mm}
\begin{abstract}
\noindent
The $R$--matrix of the symplecto--orthogonal quantum superalgebra $\Uq$ in
the vector representation is calculated, and its basic properties are
derived.
\end{abstract}

\vspace{\fill}
\noindent
math.QA/0004032

\end{titlepage}

\setcounter{page}{2}
\renewcommand{\baselinestretch}{1}
\small\normalsize

\sect{Introduction \vspace{-1ex}}
The present work is the first of two papers devoted to the construction of
the sym\-plecto--orthogonal quantum supergroup $\SP$ and of some of its
comodule superalgebras. In this work, I am going to calculate the $R$--matrix
of the sym\-plecto--orthogonal quantum superalgebra $\Uq$ in the vector
representation. Once this has been done, we can use the techniques of
Ref.~\cite{RTF} (generalized to the super case) to define the corresponding
quantum supergroup $\SP$ and to introduce its basic comodule superalgebras.
This will be carried out in the subsequent paper \cite{Sqg}.

As the reader will immediately notice, my approach is slightly different
from what he/she presumably might expect. Hence a few words of explanation
are in order. The starting point, and one of the main goals of the present
investigation, was to construct a deformed Weyl superalgebra (i.e., a
deformed oscillator algebra) $W_q(n \ty|\ty m)$, describing $n$ deformed
bosons and $m$ deformed fermions, and covariant under deformed orthosymplectic
transformations (I am grateful to V. Rittenberg for insisting that this
problem should be solved). Classically, the bosonic/fermionic commutation
relations are invariant under symplectic/orthogonal transformations. Since,
in supersymmetry, bosons/fermions are regarded to be even/odd, the natural
supersymmetric generalization of the above is that a combined system
consisting of $n$ bosons and $m$ fermions is invariant under the action of
the symplecto--orthogonal Lie superalgebra $\spo$, rather than under the
action of the orthosymplectic Lie superalgebra $\osp$. From a practical point
of view, this distinction is not really important. It is well--known that the
Lie superalgebras $\osp$ and $\spo$ are naturally isomorphic: Basically, the
transition from $\osp$ to $\spo$ amounts to a shift of the gradation of the
vector representation \cite{Sbu}. Nevertheless, I prefer to work from the
outset with the natural gradations, and to avoid any shift of gradations.

The second point where I am going to depart from the more familiar
formulation is more serious. Since Kac's basic papers on Lie superalgebras
\cite{Kcl}, \cite{Kre} it has become customary to split the family of Lie
superalgebras $\osp$ into two subfamilies, the $C$--type algebras, which are
those with $m = 1\ty$, and the $D$--type algebras, which are those with
$m \geq 2\ty\ty$. Accordingly, the so--called distinguished basis of the root
system is chosen differently for these two subfamilies.

Needless to say, there are good reasons for considering the $C$-- and
$D$--type Lie superalgebras separately. In the standard terminology, the
former are of type I\ty, while the latter are of type II\ty. This has serious
consequences for the general representation theory of these algebras. On the
other hand, one must not forget that the root systems of all of the $\osp$
algebras have some bases which resemble those of the $C$--type Lie
algebras, and others which are similar to those of the $D$--type Lie
algebras (so that I would prefer to say that these algebras are of
{\it CD}--type). In particular, for each of the $\osp$ algebras, the root system
has a basis, which is of $C$--type and contains only one odd simple root
(see Section 2). This is the basis I am going to choose (but, of course,
for the $\spo$ algebras).

Since the quantum superalgebra associated to a basic classical Lie
superalgebra depends on the choice of the basis of the root system, any such
choice has non--trivial consequences. The advantage of my choice is that it
allows of a simultaneous treatment of all cases, resulting in a unified
construction of the corresponding quantum supergroups $\SP$ and of the
deformed Weyl superalgebras $W_q(n \ty|\ty m)$. The reader might wonder
whether the differences between the $C$--type and $D$--type Lie superalgebras
will not show up at some stage of our investigations. But since in the
following we only have to consider the vector module $V$ of $\Uq$ and its
tensorial square $V \ot V$, such is not the case.

In principle, the $R$--matrix in question could be calculated by specializing
the formula for the universal $R$--matrix given in Ref.~\cite{Yam} (see also
Ref.~\cite{KTo}), or by using the results of Ref.~\cite{BSh} (I am grateful
to M. Jimbo and M. Okado for drawing my attention to the latter reference).
However, I prefer to proceed differently and to determine the corresponding
braid generator $\Rh$ by investigating the module structure of $V \ot V$.
This procedure has the advantage of yielding the spectral decomposition of 
$\Rh$ as well, moreover, at several places it can serve to check the general
theory.

The present work is set up as follows. In Section 2 we introduce the Lie
superalgebra $\spo$ and fix some notation. In particular, we specify the
basis of the root system that we are going to use, and we introduce the
corresponding Chevalley--Serre generators of the algebra. Using these data,
we define in Section 3 the quantum superalgebra $\Uq$ in the sense of
Drinfeld \cite{Dri} and Jimbo \cite{Jim} (generalized to the super case).
Basically, we follow Ref.~\cite{Yam}, but some details are different. In
Section 4 we introduce the vector module $V$ of $\Uq$. This is almost
trivial, since (in the usual sloppy terminology) this module is undeformed.
We also show that, as in the undeformed case, there exists on $V$ a
$\Uq$--invariant bilinear form, which is unique up to scalar multiples. 

In Section 5 we investigate the structure of the $\Uq$--module $V \ot V$, in
particular, we determine its module endomorphisms. This section is central
to the present work. Using the results obtained therein, we can calculate
the $R$--matrix $R$ (equivalently, the braid generator $\Rh$\ty) of $\Uq$ in
the vector representation. This will be carried out in Section 6. In
Section 7 we collect some of the basic properties of $R$ and $\Rh\ty\ty$.
Section 8 contains a comparison of our results with known special cases. A
brief discussion in Section 9 closes the main body of the paper. There are
two appendices: In Appendix A we comment on invariant bilinear forms, in
Appendix B we introduce what we have called the partial (super)transposition.

We close this introduction by explaining some of our conventions. The base
field will be the field $\CC$ of complex numbers (in the appendices, we
allow for an arbitrary field $\KK$ of characteristic zero). If $A$ is an
algebra, and if $V$ is an arbitrary (left) $A$--module, the representative
of an element $a \in A$ under the corresponding representation will be
denoted by $a_V$, and the image of an element $x \in V$ under the module
action of $a$ will be written in the form $a_V(x) = a \cd x\ty\ty$.
The multiplication in a Lie superalgebra will be denoted by pointed brackets
$\la\;\ty,\:\ra$. All algebraic notions and constructions are to be
understood in the super sense, i.e., they are assumed to be consistent with
the $\ZZ_2$\tty--\tty gradations and to include the appropriate sign factors.

\sect{Notation and a few comments on the Lie superalgebra $\spo$\vspace{-1ex}}
Essentially, we use the same type of notation as in Ref.~\cite{Sev} (see also
Ref.~\cite {Sgt}), but adapted to the present setting.

We choose two integers $m,n \geq 1$ and set
 \[ r = m + n \;. \]
Let $V = V_{\ol{0}} \op V_{\ol{1}}$ be a $\ZZ_2$\tty--\tty graded vector
space such that
 \[ \dim V_{\ol{0}} = 2n \quad,\quad \dim V_{\ol{1}} = 2m \;, \]
let $b$ be a non--degenerate, even, super--skew--symmetric, bilinear form
on $V$, and let $\spob$ be the Lie superalgebra consisting of all vector
space endomorphisms of $V$ that leave the form $b$ invariant. Then $\spob$
is isomorphic to $\spo$.

According to Ref.~\cite{Sgt}, the Lie superalgebra $\spob$ can be
described as follows (note that in the cited reference we have written
${\mathfrak{osp}}(b)$ instead of $\spob$). Let $\gl(V_{\ol{0}} \op V_{\ol{1}})$
be the general linear Lie superalgebra of the $\ZZ_2$\tty--\tty graded
vector space $V$, and let
 \[ \theta : V \ot V \lra \gl(V_{\ol{0}} \op V_{\ol{1}}) \]
be the linear map defined by
 \[ \theta(x \ot y)z = b(y\ty,z)x + (-1)^{\xi\eta} b(x\ty,z)y \;, \]
for all $x \in V_{\xi}\tty$, $y \in V_{\eta}\tty$, $z \in V$, with
$\xi,\eta \in \ZZ_2\,$. Then the kernel of $\theta$ is equal to the
subspace of all super--skew--symmetric tensors in $V \ot V$, its image is
equal to the subalgebra $\spob$ of $\gl(V_{\ol{0}} \op V_{\ol{1}})$, and
$\theta$ is an $\spob$\tty--\tty module homomorphism. In particular,
$\theta$ induces an $\spob$\tty--\tty module isomorphism of the submodule
of all super--symmetric tensors in $V \ot V$ onto the adjoint
$\spob$\tty--\tty module.

Let us make all this more explicit by introducing a suitable basis of $V$.
In order to do that we need some more notation. Define the index sets
\bea
           I \eq \{-r,-r+1,\ldots,-2\ty,-1,1,2\ty,\ldots,r-1,r\} \nn\\
  I_{\ol{0}} \eq \{-n,-n+1,\ldots,-2\ty,-1,1,2\ty,\ldots,n-1,n\} \nn\\
  I_{\ol{1}} \eq \{-r,-r+1,\ldots,-n-2\ty,-n-1,n+1,n+2\ty,\ldots,r-1,r\} \nn
\eea
and also
\[\begin{array}{rcl}
          J \eq \{-r,-r+1,\ldots,-2,-1\} \\[0.5ex]
 J_{\ol{0}} \eq \{-n,-n+1,\ldots,-2\ty,-1\} \,=\, J \cap I_{\ol{0}}
                                                                  \\[0.5ex]
 J_{\ol{1}} \eq \{-r,-r+1,\ldots,-n-2\ty,-n-1\} \,=\, J \cap I_{\ol{1}} \;.
\end{array}\]
Moreover, define the elements $\eta_i \in \ZZ_2\,$; $i \in I$, by
  \[ \eta_i = \al \quad\mbox{if $i \in I_{\al}\ty\ty$, $\al \in \ZZ_2$} \;,\]
and the sign factors
  \[ \si_i = (-1)^{\eta_i} \quad,\quad \si_{ij} = (-1)^{\eta_i\ty\eta_j}
                                        \quad\mbox{for all $i,j \in I$} \]
  \[ \tau_j = 1 \quad\mbox{and}\quad \tau_{-j} = -\si_j
                                      \quad\mbox{for all $j \in J$} \,. \]
Note that
  \[ \tau_i \ty \tau_{-i} = -\si_i \quad\mbox{for all $i \in I$} \]
(note also that the mapping $\pi : I \rar I$ used in Ref.~\cite{Sev} is given
by $\pi(i) = -i$ for all $i \in I$).

Then there exists a homogeneous basis $(e_i)_{i \in I}$ of $V$ such that
$e_i$ is homogeneous of degree $\eta_i\tty$, for all $i \in I$, and such that
 \[ b(e_i \ty, e_j) = \tau_j \ty\de_{i,-j}
                                   \quad\mbox{for all $i,j \in I$} \,. \]
We shall also use the notation
 \be C_{ij} = b(e_i \ty, e_j) \;,\;\; i,j \in I \,. \label{defC} \ee
If $C$ is the $I \ti I$--\tty matrix with elements $C_{ij}\ty\ty$, and if $G$
is the $I \ti I$--\tty matrix defined by
 \be G_{ij} = \si_i \ty \de_{ij} \quad\mbox{for all $i,j \in I$} \,,
                                                             \label{defG} \ee
we have
 \be C^2 = -G \;. \label{Csq} \ee

Besides the basis $(e_i)_{i \in I}$ of $V$, we also use the basis
$(f_i)_{i \in I}\tty$, which is dual to $(e_i)$ with respect to $b$ and is
defined by
 \[ b(f_j \ty, e_i) = \de_{ij} \quad\mbox{for all $i,j \in I$} \,. \]
Obviously, $f_i$ is homogeneous of degree $-\eta_i\tty$. Explicitly, we have
 \[ f_i \,=\, \sum_{j \in I} (C^{-1})_{ij} \ty\ty e_j
            \,=\, \tau_i \ty\ty e_{-i} \quad\mbox{for all $i \in I$} \,. \]

Using the two bases $(e_i)$ and $(f_i)$ of $V$, we define the following
elements of $\spob$:
 \[ X_{ij} = \theta(e_i \ot f_j) \quad\mbox{for all $i,j \in I$} \,. \]
Let $(E_{ij})_{i,\ty j \in I}$ be the basis of
$\gl(V_{\ol{0}} \op V_{\ol{1}})$ that canonically corresponds to the basis
$(e_i)_{i \in I}$ of $V$, i.e.,
 \[ E_{ij}(e_k) = \de_{jk} \ty\ty e_i
                                    \quad\mbox{for all $i,j,k \in I$} \,.\] 
Then we obtain
 \[ X_{ij}
 \,=\, E_{ij} + \si_{ij}\sum_{k,\ty\ell}
                                 C_{ik}\ty(C^{-1})_{j\ell}\tty\ty E_{\ell k}
 \,=\, E_{ij} + \si_{ij}\ty\tau_{-i}\ty\tau_{j}\ty E_{-j,-i}
                                       \quad\mbox{for all $i,j \in I$} \,. \]
In particular, we have
 \[ X_{ii} = E_{ii} - E_{-i,-i} \quad\mbox{for all $i \in I$} \,. \]
According to the properties of the map $\theta\ty$, the elements $X_{ij}$
generate the vector space $\spob$, moreover, the super--symmetry of
$\theta$ implies that
 \[ \tau_{-j} \tty X_{i,-j} \,=\, \si_{ij} \ty \tau_{-i} \tty X_{j,-i}
                                      \quad\mbox{for all $i,j \in I$} \,. \]
Thus we have
 \[ X_{i,-i} = 0 \quad\mbox{for all $i \in I_{\ol{1}}$} \;. \]

Let $\fh$ be the subspace of $\spob$ that is spanned by the elements
$X_{ii}\ty\ty$, $i \in I$. Obviously, $\fh$ consists of those elements of
$\spob$ whose matrices with respect to the basis $(e_i)$ are diagonal, and
the $X_{jj}$ with $j \in J$ form a basis of $\fh\ty$.

Define, for every $i \in I$, the linear form $\ve_i$ on $\fh$ by
 \[ H(e_i) = \ve_i(H) \ty e_i \quad\mbox{for all $H \in \fh$} \;. \]
Then it is easy to see that
 \[ \ve_{-i} = - \ve_i \quad\mbox{for all $i \in I$} \,, \]
and that
 \[ \ve_i(X_{jj}) = \de_{ij} \quad\mbox{for all $i,j \in J$} \,. \]
Thus $(\ve_j)_{j \in J}$ is the basis of $\fh^{\ast}$ that is dual to the
basis $(X_{jj})_{j \in J}$ of $\fh\ty$.

Since $\theta$ is an $\spob$\tty--\tty module homomorphism, it follows that
 \[ \la H, X_{ij} \ra \,=\, (\ve_i - \ve_j)(H) \ty X_{ij} \]
for all $H \in \fh$ and all $i,j \in I$ (recall that the multiplication in a
Lie superalgebra is denoted by pointed brackets). We conclude that $\fh$ is
a Cartan subalgebra of $\spob$, that
  \be \De = \{\ve_i - \ve_j \ty|\, i,j \in I;\;
                     \mbox{$j \neq i$ and $j < -i \ty,\ty$
                            or $j = -i \in I_{\ol{0}}\ty$}\} \label{rs} \ee
is the root system of $\spob$ with respect to $\fh\ty$, and that $X_{ij}$ is
a (non--zero) root vector corresponding to the root $\ve_i - \ve_j$ (with
$i,j$ as specified on the right hand side of Eqn.~\reff{rs}). The root
$\ve_i - \ve_j$ is even/odd depending on whether
$\si_i \ty\ty \si_j = \pm 1\ty$.

In order to introduce an adequate bilinear form on $\fh^{\ast}$, we recall
that the invariant bilinear form
 \[ (X,Y) \lra \ohal \tty \mbox{Str}(X Y) \]
on $\spob$ is non--degenerate and super--symmetric; consequently, its
restriction to $\fh$ is likewise. Let $(\;\,|\;\,)$ denote the bilinear form
on $\fh^{\ast}$ that is inverse to this restriction. By definition, we have
 \[ (\lam \ty|\ty \mu) \,=\, \ohal \tty \mbox{Str}(H_{\lam} \ty H_{\mu}) \]
for all $\lam\ty,\mu \in \fh^{\ast}$, where, for example, the element
$H_{\lam} \in \fh$ is uniquely determined through the equation
 \be \lam(H) \,=\, \ohal \tty \mbox{Str}(H_{\lam} \ty H)
                      \quad\mbox{for all $H \in \fh$} \;.  \label{Hlam} \ee
It is easy to check that
 \[ H_{\lam} \,=\, \sum_{j \in J} \si_j \tty \lam(X_{jj}) \ty X_{jj}
                           \quad\mbox{for all $\lam \in \fh^{\ast}$} \,, \]
and that
 \[ (\ve_i \ty|\ty \ve_j) \,=\, \si_i \tty \de_{ij}
                                     \quad\mbox{for all $i,j \in J$} \,. \]

Let us now specify the basis of the root system $\De$ that we are going to
use in the following. It is equal to $(\al_j)_{j \in J}\tty$, where the
simple roots $\al_j$ are defined by
  \[  \al_j = \left\{ \begin{array}{ll}
                     \ve_j - \ve_{j+1} & \mbox{for $-r \leq j \leq -2$} \\
                            2 \ve_{-1} & \mbox{for $j = -1$} \;.
                      \end{array} \right. \]
Note that $\al_{-n-1}$ is the sole odd simple root. The corresponding
Chevalley--Serre generators of $\spob$ are denoted by $E_j^v \ty$,
$F_j^v \ty$, $H_j^v \ty\ty$; $j \in J$, and are introduced as follows.
First of all, we choose
  \[ \begin{array}{rcl}
  E_j^v \eq \left\{ \begin{array}{ccll}
                     X_{j,\ty j+1}
                         \eqq E_{j,\ty j+1}
                              - \si_j\ty\ty\si_{j,\ty j+1}\ty E_{-j-1,-j}
                               & \hspace{0.95em}
                                 \mbox{for $-r \leq j \leq -2$} \\[0.5ex]
                     \ohal X_{-1,\ty 1} \eqq E_{-1,\ty 1}
                               & \hspace{0.95em}
                                 \mbox{for $j = -1$}
                  \end{array} \right. \vspace{1ex} \\
  F_j^v \eq \left\{ \begin{array}{ccll}
                     X_{j+1,\ty j}
                         \eqq E_{j+1,\ty j}
                             - \si_{j+1}\ty\si_{j+1,\ty j}\ty E_{-j,-j-1}
                               & \mbox{for $-r \leq j \leq -2$} \\[0.5ex]
                     \ohal X_{1,-1} \eqq E_{1,-1}
                               & \mbox{for $j = -1$} \;.
                    \end{array} \right.
     \end{array} \vspace{2.0ex} \]

\noindent
{\it Remark 2.1.} It is easy to check that
 \be \si_{j,\ty j+1} = \si_{j+1} \quad\mbox{for $-r \leq j \leq -2$} \;,
                                                       \label{relsione} \ee
or, equivalently, that
 \be \si_{j,\ty j+1} = \si_j \quad\mbox{for $1 \leq j \leq r-1$} \,.
                                                       \label{relsitwo} \ee
This implies that in the equation for $E_j^v \ty$, $-r \leq j \leq -2\ty\ty$,
the factor $\si_{j,\ty j+1}$ might be replaced by $\si_{j+1}\ty$, and in
the equation for $F_j^v \ty$, $-r \leq j \leq -2\ty\ty$, the $\si$--factors
might be dropped. I prefer to keep the $\si$--factors as they stand: They
have an immediate meaning in terms of the sign rules of supersymmetry, and
by modifying the equations for the $E_j^v$ and $F_j^v$ as mentioned
above, we might well end up in the unpleasant situation where we would
have to check (possibly implicitly) the equations \reff{relsione} or
\reff{relsitwo} again and again. \vspace{2.0ex}

Using the elements $E_j^v$ and $F_j^v\ty$, we define the generators
$H_j^v \in \fh$ as usual by
  \be  H_j^v = \la E_j^v , F_j^v \ra \quad\mbox{for all $j \in J$} \,.
                                                            \label{Hjv} \ee
More explicitly, we find
  \[  H_j^v \,=\, \left\{ \begin{array}{ll}
                  X_{jj} - \si_j \ty \si_{j+1} \ty X_{j+1,\ty j+1}
                                & \mbox{for $-r \leq j \leq -2$} \\[0.5ex]
                      X_{-1,-1} & \mbox{for $j = -1$} \;.
                        \end{array} \right. \]
Then the generators $E_j^v \ty$, $F_j^v \ty$, $H_j^v$ satisfy the following
familiar relations, which hold for all $i,j \in J \,$:
\be   \la H_i^v \ty, H_j^v \ra = 0  \label{hh}  \ee
\be   \la E_i^v \ty, F_j^v \ra = \de_{ij} \ty H_j^v  \label{efr}  \ee
\be   \la H_i^v \ty, E_j^v \ra = a_{ij}\ty E_j^v \quad,\quad
            \la H_i^v \ty, F_j^v \ra = -a_{ij}\ty F_j^v \,. \label{hef} \ee
Here, $A = (a_{ij})_{i,\ty j \in J}$ is the Cartan matrix, whose elements are
defined by
  \[  a_{ij} = \al_j(H_i^v) \quad\mbox{for all $i,j \in J$} \,. \]
The Cartan matrix is tridiagonal. For $n \geq 2\tty$, it takes the following
form:
  \[ A \;=\; \left( \begin{array}{rrrrrrrrrrrr}
    2 & -1 &    & \mc{9}{c}{} \\
   -1 &  2 & -1 & \mc{9}{c}{} \\
   \mc{12}{c}{} \\
   \mc{12}{c}{} \\
   \mc{3}{c}{} & -1 &      2 &     -1 & \mc{6}{c}{} \\
   \mc{4}{c}{} & -1 &      0 & \;\; 1 & \mc{5}{c}{} \\
   \mc{5}{c}{} & -1 & \;\; 2 &     -1 & \mc{4}{c}{} \\
   \mc{12}{c}{} \\
   \mc{12}{c}{} \\
   \mc{8}{c}{} & -1 &  2 & -1 &    \\
   \mc{8}{c}{} &    & -1 &  2 & -2 \\
   \mc{8}{c}{} &    &    & -1 &  2 
                    \end{array} \right) \,, \]
where the zero on the diagonal has the row and column number $-n-1\ty$. For
$n = 1\ty$, the zero is in position $(-2\tty,-2)$, and the lower right
corner of $A$ is equal to
  \[ \left( \begin{array}{rrr}
                       2 & -1 & \;\; 0 \\
                      -1 &  0 & \;\; 2 \\
                       0 & -1 & \;\; 2 
            \end{array} \right) \,. \] 
Finally, for $m = n = 1$ the Cartan matrix is given by
  \[ A \,=\, \left( \begin{array}{rr}
                            0 & 2 \\
                           -1 & 2
                    \end{array} \right) \,. \]
Thus the Dynkin diagram  of $\spob \simeq \spo$ with respect to our
basis of the root system takes the form
\vspace{-2.0ex}

\begin{center}
\unitlength1.25mm
\begin{picture}(116,10)
\thicklines
\put(2.2,5){\circle{4.4}}
\put(0.0,-0.5){$-r$}
\put(4.4,5){\line(1,0){10}}
\put(16.6,5){\circle{4.4}}
\put(11.2,-0.5){$-r+1$}
\put(18.8,5){\line(1,0){10}}
\multiput(30.3,4.8)(2.0,0){4}{.}
\put(38.8,5){\line(1,0){10}}
\put(51.0,5){\circle{4.4}}
\put(45.4,-0.5){$-n-1$}
\put(49.6,3.5){\line(1,1){2.9}}
\put(49.6,6.5){\line(1,-1){2.9}}
\put(53.2,5){\line(1,0){10}}
\put(65.4,5){\circle{4.4}}
\put(63.1,-0.5){$-n$}
\put(67.6,5){\line(1,0){10}}
\multiput(79.1,4.8)(2.0,0){4}{.}
\put(87.6,5){\line(1,0){10}}
\put(99.8,5){\circle{4.4}}
\put(97.5,-0.5){$-2$}
\put(102.0,5){\line(1,1){2.9}}
\put(102.0,5){\line(1,-1){2.9}}
\put(102.6,5.5){\line(1,0){9.4}}
\put(102.6,4.5){\line(1,0){9.4}}
\put(114.2,5){\circle{4.4}}
\put(112.2,-0.5){$-1$}
\end{picture}
\end{center}
\vspace{0.1ex}
\centerline{Dynkin diagram of the Lie superalgebra $\spo$}
\vspace{2.5ex}

\noindent
{\it Remark 2.2.} It may be helpful to comment on the rules according to which
the generators $H_i^v$ (and hence $E_i^v$ and $F_i^v$) have been chosen. If
the simple root $\al_i$ is even, we choose $H_i^v$ such that
$\al_i(H_i^v) = 2\ty\ty$. For odd simple roots, the situation is more
complicated. In the present case, the sole odd simple root $\al_i\:$,
$i = -n-1\ty$, is such that $(\al_i \ty|\ty \al_i) = 0\ty$. Then it follows
that $a_{ii} = 0\ty$, and the element $H_i^v$ is usually chosen such that
(for this index $i$)
 \[ a_{ij} \in \ZZ \quad\mbox{for all $j \in J$} \,, \]
and such that these $a_{ij}$ don't have a common divisor. This fixes the
$a_{ij}$ up to a common sign factor, which (according to Kac) is chosen such
that $a_{i,\ty i+1} > 0$ (assuming that $i+1 \in J$ and that
$a_{i,\ty i+1} \neq 0$). These conventions are introduced simply for
convenience, and they are of little (if any) importance. Note that, for
$m = n = 1\ty$, we haven't followed these conventions: The first row of the
Cartan matrix could be divided by $2\tty$, and for
$\sli(2 \ty|\ty 1) \simeq \spott$ this is usually done. Our choice is
motivated by the wish for a unified treatment of all cases. \vspace{2.0ex}

The relations \reff{hh} -- \reff{hef} given above are not sufficient to
characterize the Lie superalgebra $\spo$ completely, there are certain
Serre--type and supplementary relations which must also be satisfied. We
don't give these relations here, but only mention that they can be read
off from the relations \reff{gkom} -- \reff{zus7} by setting $q = 1\ty$.

\sect{Definition of the quantum superalgebra \\ $\Uq$ \vspace{-1ex}}
The notation introduced in the preceding section will now be used to define
the quantum superalgebra $\Uq$. Basically, we are going to follow
Ref.~\cite{Yam}, however, there will be differences in detail.

Define the diagonal $J \ti J$--\tty matrix $D$ by
 \[ D \,=\, (d_i \ty\ty \de_{ij})_{i,\ty j \in J} \,=\,
     \mbox{diag}(\ub{-1,-1,\ldots,-1}_m , \ub{1,1,\ldots,1}_{n-1} , 2) \;. \]
It is chosen such that
 \[ (DA)_{ij} = (\al_i \ty|\ty \al_j) \quad\mbox{for all $i,j \in J$} \,.\]
In particular, the matrix $DA$ is symmetric.

Let $q \in \CC$ be a non--zero complex number, and assume that $q$ is
{\em not a root of unity.} We use the abbreviation
 \[ q_i = q^{d_i} \;. \]
Then we have
 \[ q_i^{a_{ij}} = q^{(\al_i |\ty \al_j)}
                                      \quad\mbox{for all $i,j \in J$} \,.\]
Now the quantum superalgebra $\Uq$ is defined to be the universal associative
superalgebra (with unit element) with generators $K_i$\,, $K_i^{-1}$,
$E_i$\,, $F_i$\,; $i \in J$, and certain relations to be specified below. The
$\ZZ_2$\tty--\tty gradation is fixed by requiring that $E_{-n-1}$ and
$F_{-n-1}$ be odd, while all the other generators are even. (Needless to
say, one has to check that the relations are compatible with this
requirement.) The relations are the following, they are assumed to hold for
all $i,j \in J$\,:
  \be K_i \ty K_i^{-1} = K_i^{-1} K_i = 1 \vspace{0.5ex} \label{inv} \ee
  \be K_i \ty K_j = K_j \ty K_i \vspace{0.5ex} \label{kkom} \ee
  \be K_i \ty E_j \ty K_i^{-1} = q_i^{a_{ij}} E_j \quad,\quad
     K_i \ty F_j \ty K_i^{-1} = q_i^{-a_{ij}} F_j \vspace{0.5ex}
                                                             \label{kef} \ee
  \be \la E_i \ty, F_j \ra = \de_{ij}\,\frac{K_i - K_i^{-1}}{q_i - q_i^{-1}}
                                          \;.  \vspace{0.5ex} \label{ef} \ee
In addition, the generators $E_i$ satisfy certain Serre--type and
supplementary relations among themselves, as do the generators $F_i\ty\ty$.
We only write the relations for the $E_i\ty\ty$, those for the $F_i$ are
obtained from these by simply replacing $E$ by $F$.

In the subsequent relations, it is always assumed that $i,j \in J$. Suppose
first that the root $\al_i$ is even, i.e., that $i \neq -n-1$\,. Then we have
  \be  \la E_i \ty, E_j \ra = 0 \quad\mbox{for $a_{ij} = 0$} \label{gkom} \ee
  \be  E_i^2 \ty E_j - (q + q^{-1})\ty E_i \ty E_j \ty E_i + E_j \ty E_i^2 = 0
        \quad\mbox{for $i \leq -3$\,, $|i - j| = 1$} \,. \label{serij} \ee
If $\al_{-2}$ is even, i.e., if $n \geq 2$\,, we have (as in the case of
symplectic Lie algebras) 
  \be  E_{-2}^2 \ty E_{-3} - (q + q^{-1})\ty E_{-2} \ty E_{-3} \ty E_{-2} 
                                      + E_{-3} E_{-2}^2 = 0 \label{ser23} \ee
  \be  E_{-2}^3 \ty E_{-1} 
            - (q^2 + 1 + q^{-2})\ty E_{-2}^2 \ty E_{-1} \ty E_{-2}
            + (q^2 + 1 + q^{-2})\ty E_{-2} \ty E_{-1} \ty E_{-2}^2
                              - E_{-1} \ty E_{-2}^3 = 0 \,. \label{ser21} \ee 
In all cases, the generators $E_{-2}$ and $E_{-1}$ satisfy
  \be  E_{-1}^2 \ty E_{-2} - (q^2 + q^{-2})\ty E_{-1} \ty E_{-2} \ty E_{-1}
                              + E_{-2} \ty E_{-1}^2 = 0 \,. \label{ser12} \ee
Next we recall that $\al_{-n-1}$ is the sole odd simple root, and that this
root is isotropic. Correspondingly, we have
  \be  \la E_{-n-1} \ty, E_j \ra = 0 \quad\mbox{for $a_{-n-1,\ty j} = 0$}
                                                        \,, \label{ugkom} \ee
in particular,
  \be  E_{-n-1}^2 = 0 \,. \label{nilq} \ee
Finally, there are the following supplementary relations. If $m,n \geq 2$\,,
we have
\be  \la E_{-n-1} \ty,\la E_{-n} \ty,\la E_{-n-1} \ty, E_{-n-2}
                             \ra_q\ty\ra_{q^{-1}}\ra = 0 \,, \label{zus4} \ee
and for $n = 1\ty$, $m \geq 3$ we have
\be  \la E_{-2} \ty,\la E_{-3} \ty,\la E_{-2} \ty,\la E_{-1} \ty,\la
                     E_{-2} \ty,\la E_{-3} \ty, E_{-4}
 \ra_q\ty\ra_q\ty\ra_{q^{-2}}\ra_{q^{-1}}\ra_q\ty\ra = 0 \,. \label{zus7} \ee
The last two relations are expressed in terms of so--called
$q$--supercommutators. We recall the definition: If $A$ is any associative
superalgebra, if $p$ is any non--zero complex number, and if $X \in A_{\xi}$
and $Y \in A_{\eta}$\,, with $\xi,\eta \in \ZZ_2\,$, the
$p$--supercommutator of $X$ and $Y$ is defined by
 \[ \la X \ty, Y \ty\ra_p
                     = X \ty Y - p \tty (-1)^{\xi \ty \eta}\ty Y \ty X \;.\]
Obviously, we have
 \[ \la X \ty, Y \ty\ra_1 = \la X \ty, Y \ty\ra \;. \]
As shown below, the Serre--type relations can also be expressed in terms of
$q$--supercommutators.

The superalgebra $\Uq$ is converted into a Hopf superalgebra by means of
structure maps, which are fixed by the following equations: \\[0.5ex]
coproduct
  \[ \begin{array}{rcl}
        \De(K_i^{\pm 1}) \eqq K_i^{\pm 1} \ot K_i^{\pm 1} \\[1.0ex]
                \De(E_i) \eqq E_i \ot 1 + K_i \ot E_i \\[1.0ex]
                \De(F_i) \eqq F_i \ot K_i^{-1} + 1 \ot F_i 
     \end{array} \]
counit
  \[ \ve(K_i^{\pm 1}) = 1 \quad,\quad \ve(E_i) = \ve(F_i) = 0 \]
antipode
  \[ S(K_i^{\pm 1}) = K_i^{\mp 1} \quad,\quad
            S(E_i) = - K_i^{-1} E_i \quad,\quad S(F_i) = -F_i \ty K_i \;.\]

In the subsequent series of remarks, we collect some elementary properties
of the Hopf superalgebra $\Uq$. \vspace{2.0ex}

\noindent
{\em Remark 3.1.} Let $\bf Q$ be the root lattice of $\spo$, i.e.,
  \[ {\bf Q} = \sum_{i \in J} \ZZ\tty\al_i \;. \]
Then the algebra $\Uq$ admits a unique $\bf Q$\tty--\tty gradation such that,
for all $i \in J$, the element $E_i$ is homogeneous of degree $\al_i\ty\ty$,
$F_i$ is homogeneous of degree $-\al_i\ty\ty$, and $K_i$ is homogeneous of
degree $0\tty$. In view of a more general definition to be given later, the
$\bf Q$\tty--\tty degree of a $\bf Q$\tty--\tty homogeneous element is called
its {\em weight.} If an element $X \in \Uq$ has the weight
$\lam \in {\bf Q}\ty\ty$, it satisfies
  \[ K_j \ty X \ty K_j^{-1} = q^{(\al_j |\ty \lam)} \ty X
                                       \quad\mbox{for all $j \in J$} \,. \]
Conversely, if an element $X \in \Uq$ satisfies this condition, it is
$\bf Q$\tty--\tty homogeneous of weight $\lam$ (since $q$ is not a root of
unity). Note that the structure maps $\De\,$, $\ve\ty\ty$, and $S$ are
$\bf Q$\tty--\tty homogeneous of degree zero. \vspace{2.0ex}

\noindent
{\em Remark 3.2.} The antipode $S$ is bijective. To prove this, we show that
$S^2$ is bijective. Indeed, it is easy to check that, for all $i \in J$,
 \be S^2(K_i) = K_i
        \quad,\quad S^2(E_i) = q^{-(\al_i |\ty \al_i)} \ty E_i
              \quad,\quad S^2(F_i) = q^{(\al_i |\ty \al_i)} \ty F_i \;.
                                                              \label{sq} \ee
Since $S^2$ is an algebra endomorphism of $\Uq$, and since a suitable set of
monomials in the generators $K_i^{\pm 1}$, $E_i\ty\ty$, and $F_i$ forms a
basis of $\Uq$, this implies our claim.

Actually, $S^2$ is an inner automorphism of the algebra $\Uq$. Let $2\ty\rho$
denote the sum of the even positive roots minus the sum of the odd positive
roots of $\spo$. Explicitly, we have
  \[ 2\ty\rho = -2 \!\sum_{i = -r}^{-n-1}(i+2n+1)\ty\ve_i
                -2 \!\sum_{i = -n}^{-1} i\tty\ve_i \;, \]
and it is easy to check that
  \[ (2\ty\rho \ty|\ty \al_i) = (\al_i \ty|\ty \al_i)
                                    \quad\mbox{for all $i \in J$} \,. \]
Given an arbitrary linear combination of the simple roots $\al_i$ with
coefficients $r_i \in \ZZ$\,,
  \[ \lam = \sum_{i \in J} r_i\ty\al_i \in {\bf Q} \;, \]
we define
  \[ K_{\lam} = \prod_{i \in J} K_i^{r_i} \,. \]
In particular, we have
  \[ K_{\al_i} = K_i \quad\mbox{for all $i \in J$} \,. \]
Then the Eqns.~\reff{sq} immediately imply that 
  \[ S^2(X) = K_{-2\ty\rho} \ty\ty X \ty K_{-2\ty\rho}^{-1}
                       \quad\mbox{for all $X \in \Uq$} \;. \vspace{0.5ex} \]

\noindent
{\em Remark 3.3.} Obviously, there is a certain symmetry between the $E$ and
$F$ generators of $\Uq$. To make this more explicit, we note that there is a
unique algebra endomorphismus
  \[ \vp : \Uq \lra \Uq \;, \]
such that for all $i \in J$
  \[ \vp(E_i) = F_i \quad,\quad \vp(F_i) = (-1)^{\ga_i} E_i
                           \quad,\quad \vp(K_i^{\pm 1}) = K_i^{\mp 1} \;, \]
where $\ga_i \in \ZZ_2$ is the degree of $E_i\ty\ty$. This endomorphism is
homogeneous of $\ZZ_2$\tty--\tty degree zero, and we have $\vp^4 = \id$\,.
Consequently, $\vp$ is an automorphism of the (associative) superalgebra
$\Uq$.

Now let $\Uq^{\mbox{\scriptsize cop}}$ be the bi--superalgebra which,
regarded as a $\ZZ_2$\tty--\tty graded algebra, coincides with $\Uq$, but
whose coalgebra structure is opposite (in the super sense) to that of $\Uq$.
Then it is easy to check that
  \be \vp : \Uq \lra \Uq^{\mbox{\scriptsize cop}} \label{phi} \ee
is a homomorphism of bi--superalgebras. Since $\vp$ is bijective, this
implies that $\Uq^{\mbox{\scriptsize cop}}$ is a Hopf superalgebra, and that
$\vp$ is a Hopf superalgebra iso\-morphism. As is well-known, it follows
(once again) that $S$ is bijective, and that $S^{-1}$ is the antipode of
$\Uq^{\mbox{\scriptsize cop}}$. \vspace{2ex}

The Serre--type and the supplementary relations can be written in various
ways. Before we do that, we remind the reader of the definition of the
adjoint representation of a Hopf superalgebra $H$. This is a (graded)
representation of the superalgebra $H$ on the graded vector space $H$, it
is denoted by ${\ad}\tty$, and is defined as follows. Let $X$ be an
arbitrary element of $H$, and set
 \[ \De(X) = \sum_\au X_a^1 \ot X_a^2 \;, \vso \]
with homogeneous elements $X_a^1 \ty, X_a^2 \in H$, of degree $\xi_a^1$
and $\xi_a^2 \ty$, respectively. Then $\add X$ (the representative of $X$)
is given by
  \[ (\add X)(Y)
        \,=\, \sum_\au (-1)^{\xi^2_a \ty \eta}\ty X^1_a \ty Y S(X^2_a) \;,
                                                                     \vso \]
for all homogeneous elements $Y \in H_{\eta}\ty\ty$, where
$\eta \in \ZZ_2\,$. We note that $\add X$ is a generalized derivation in
the sense that, if $Y'$ is another element of $H$, 
  \[ (\add X)(Y \ty Y')
           \,=\, \sum_\au (-1)^{\xi^2_a \ty \eta}
                           (\add X^1_a)(Y) \ty (\add X^2_a)(Y') \;. \vso \] 
Now suppose that $S$ bijective. This implies that $H^{\mbox{\scriptsize cop}}$
(see the analogous definition of $\Uq^{\mbox{\scriptsize cop}}$ given above)
is a Hopf superalgebra with antipode $S^{-1}$. Let $\adb$ be the adjoint
representation of $H^{\mbox{\scriptsize cop}}$. Then $\adb$ is a graded
representation of $H$ in $H$, it is given by 
  \[ (\addb X)(Y)
          \,=\, \sum_\au (-1)^{\xi^1_a (\xi^2_a + \eta)}\ty
                                      X^2_a \tty Y S^{-1}(X^1_a) \;,\vso \]
and it satisfies
  \[ (\addb X)(Y \ty Y')
             \,=\, \sum_\au (-1)^{\xi^1_a (\xi^2_a + \eta)}
                         (\addb X^2_a)(Y) \ty (\addb X^1_a)(Y') \;. \vso \] 
Let us now choose $H = \Uq$. Then the isomorphism $\vp$ given in
\reff{phi} shows that
 \be \addb(\vp(X)) = \vp \comp (\add X) \comp \vp^{-1}
                       \quad\mbox{for all $X \in \Uq$} \;. \label{adphi} \ee
Moreover, for every element $X \in \Uq$ of weight $\lam\ty\ty$ we have
 \bea (\add E_i)(X) \eqq \la E_i \ty, X \ra_{q^{(\al_i |\ty \lam)}}
                                                       \label{adE} \\[0.5ex]
     (\addb F_i)(X) \eqq \la F_i \ty, X \ra_{q^{-(\al_i |\ty \lam)}} \;.
                                                           \label{adbF} \eea
We note that in the proof of these equations we only have to use the first
resp.~second of the relations \reff{kef} but none of the other defining
relations.

Now Eqn.~\reff{adE} implies that the left hand side of Eqn.~\reff{gkom} is
equal to
  \[ \la E_i \ty, E_j \ra = (\add E_i)(E_j) \;, \]
the left hand side of Eqn.~\reff{serij} is equal to
  \[ \la E_i \ty,\la E_i \ty, E_j \ra_{q^{\mp 1}}\ra_{q^{\pm 1}}
                                               = (\add E_i)^2\ty(E_j) \;, \]
the left hand side of Eqn.~\reff{ser23} is equal to
  \[ \la E_{-2} \ty,\la E_{-2} \ty, E_{-3} \ra_{q^{\mp 1}}\ra_{q^{\pm 1}}
                                         = (\add E_{-2})^2\ty(E_{-3}) \;, \]
the left hand side of Eqn.~\reff{ser21} is equal to
  \[ \la E_{-2} \ty,\la E_{-2} \ty,\la E_{-2} \ty, E_{-1}
                              \ra_{q^{\mp 2}}\ty\ra\ty\ra_{q^{\pm 2}}
                                         = (\add E_{-2})^3\ty(E_{-1}) \;, \]
the left hand side of Eqn.~\reff{ser12} is equal to
  \[ \la E_{-1} \ty,\la E_{-1} \ty, E_{-2} \ra_{q^{\mp 2}}\ra_{q^{\pm 2}}
                                         = (\add E_{-1})^2\ty(E_{-2}) \;, \]
the left hand side of Eqn.~\reff{ugkom} is equal to
  \[ \la E_{-n-1} \ty, E_j \ra = (\add E_{-n-1})(E_j) \;, \]
the left hand side of Eqn.~\reff{zus4} is equal to
\begin{eqnarray}
\lefteqn{\hspace*{-2em}
   \la E_{-n-1} \ty,\la E_{-n} \ty,\la E_{-n-1} \ty, E_{-n-2}
                                                  \ra_q\ty\ra_{q^{-1}}\ra
                                             \hspace{10.5em}} \nn \\[0.5ex]
   & \hspace{10.5em} = &\!\! 
                 (\add E_{-n-1})(\add E_{-n})(\add E_{-n-1})(E_{-n-2}) \;,
                                             \label{zus4var} \end{eqnarray}
and the left hand side of Eqn.~\reff{zus7} is equal to
\begin{eqnarray}
\lefteqn{\hspace*{-2em}
   \la E_{-2} \ty,\la E_{-3} \ty,\la E_{-2} \ty,\la E_{-1} \ty,\la
                                     E_{-2} \ty,\la E_{-3} \ty, E_{-4}
                      \ra_q\ty\ra_q\ty\ra_{q^{-2}}\ra_{q^{-1}}\ra_q\ty\ra
                                              \hspace{2.9em}} \nn \\[0.5ex]
   & \hspace{2.9em} = &\!\! (\add E_{-2})(\add E_{-3})(\add E_{-2})
                        (\add E_{-1})(\add E_{-2})(\add E_{-3})(E_{-4}) \;.
                          \st{-1.0ex}{1.0ex} \label{zus7var} \end{eqnarray}

\noindent
{\em Remark 3.4.} Using Eqn.~\reff{nilq} and the fact that $E_{-n-2}$ and
$E_{-n}$ commute, it is easy to see that the left hand side of
Eqn.~\reff{zus4var} is invariant under the substitution $q \rar q^{-1}$.
Somewhat unexpectedly, it seems that this is not the case for the left hand
side of Eqn.~\reff{zus7var}, even if one assumes that all the relations for
the $E$--generators except Eqn.~\reff{zus7} are satisfied. \vspace{2.0ex}

\noindent
{\em Remark 3.5.} Taking the defining relations for granted except
\reff{zus7}, one can show that the expressions in Eqn.~\reff{zus7var} are
annihilated by all $\add F_i\ty\ty$. Since, quite generally, we have
  \[ (\add F_i)(X) = \la F_i \ty, X \ra \ty K_i 
                   \quad\mbox{for all $i \in J$ and all $X \in \Uq$} \;, \]
the same is true when acting with $\la F_i \ty, \ty\cd\, \ra$. This
shows that by ``acting'' on the relation \reff{zus7} with the generators
$F_i\ty\ty$, we cannot derive new relations for the \vspace{2.0ex}
$E$--generators.

Up to now we have only discussed the Serre--type and supplementary relations 
for the $E$--generators. Of course, similar comments can be made for the
$F$--generators as well, but with $\ad$ replaced by $\addb$. In fact, all we 
have to do is to apply the isomorphism $\vp$ given in \reff{phi} and to
recall Eqn.~\reff{adphi}.

We close this section by a remark on the weights of a $\Uq$--module $W$.
In the present work, all $\Uq$--modules will be weight modules, in the sense
that the representatives $(K_j)_W\ty\ty$, $j \in J$, are simultaneously
diagonalizable, and such that, for any common eigenvector $x$ of these
operators, we have
  \[ K_j \cd x = q^{(\al_j |\ty \lam)} \ty x 
                                       \quad\mbox{for all $j \in J$} \,, \]
with a linear form
  \[ \lam \in \sum_{i \in J} \ZZ\tty\ve_i \;. \]
Since $q$ is not a root of unity, the linear form $\lam$ is uniquely fixed
by these conditions and is called the {\em weight} of $x\ty\ty$. This
definition generalizes the definition of the weight of an element of $\Uq$:
In that case the representation considered is the adjoint representation
$\ad$ or its modified version $\addb$.

\sect{The vector module $V$ of $\Uq$ \vspace{-1ex}}
Let us now discuss the vector module $V$ of $\Uq$. The definition of $V$ is
easy, since (in the usual sloppy terminology) the vector module of $\Uq$ is
undeformed. More precisely, let $V$ be the graded vector space introduced in
Section 2, and let $E_i^v \ty$, $F_i^v \ty$, and $H_i^v \ty\ty$; $i \in J$,
be the linear operators on $V$ defined there. Define the linear operators
$K_i^v\ty\ty$, $i \in J$, by
 \[ K_i^v = q_i^{H_i^v} \quad\mbox{for all $i \in J$} \,. \]
Since the operators $H_i^v$ are diagonalizable, with eigenvalues
$0, \pm 1\ty$, the operators $K_i^v$ are well--defined and, obviously,
invertible. It is easy to see that the operators $E_i^v \ty$, $F_i^v \ty$,
and $(K_i^v)^{\pm 1}$ satisfy the defining relations of the generators
$E_i$\,, $F_i$\,, and $K_i^{\pm 1}\ty\ty$; $i \in J$. Hence there exists
a unique graded representation $\pi$ of the algebra $\Uq$ in $V$ such that
 \[ \pi(E_i) = E_i^v \quad,\quad
                \pi(F_i) = F_i^v \quad,\quad
                            \pi(K_i^{\pm 1}) = (K_i^v)^{\pm 1}
                                       \quad\mbox{for all $i \in J$} \,. \]
The graded vector space $V$, endowed with this representation, will be
called the vector module of $\Uq$. \vspace{2.0ex}

\noindent
{\em Remark 4.1.} The reader might suspect that checking the seventh order
relation \reff{zus7} might be quite tedious. Actually, this is not the case.
Let $\Lgr(V)$ be the superalgebra of all linear operators in $V$. It is
well--known that $\Lgr(V)$ is an $\spob$--module in a canonical way, and its
weights (with respect to the Cartan subalgebra $\fh$) are the linear forms
$\ve_i - \ve_j$\,; $i,j \in I$. Since $\spob$ acts on $\Lgr(V)$ by
superderivations, any product with one factor $E_{-1}^v\ty$, three factors
$E_{-2}^v\ty$, two factors $E_{-3}^v\ty$, and one factor $E_{-4}^v$ has
the weight $\ve_{-4} + \ve_{-3} + \ve_{-2} - \ve_{-1}$\,. Since this is not
a weight of $\Lgr(V)$, every such product is equal to zero, and this
implies the relation to be proved. \vspace{2.0ex}

Let $(e_i)_{i \in I}$ be the basis of $V$ used in Section 2. Then we have
  \[ K_j \cd e_i = q^{(\al_j |\ty \ve_i)} \ty e_i
                     \quad\mbox{for all $j \in J$ and all $i \in I$} \,. \]
Stated differently, $e_i$ is a weight vector with weight $\ve_i\ty\ty$,
just as in the undeformed case.

Our next goal is to show that there exists a unique (up to scalar multiples)
$\Uq$--invariant bilinear form on $V$. (For a few comments on invariant
bilinear forms, see Appendix A.) Let $b$ be a bilinear form on $V$, and
let $\bt$ be the linear form on $V \ot V$ canonically corresponding to
$b\ty\ty$. Then $b$ is $\Uq$--invariant if and only if
  \[ \bt\tty(X \cd (x \ot y)) = \ve(X)\ty\ty \bt\tty(x \ot y) \]
for all $X\in \Uq$ and all $x,y \in V$ (see Eqn.~\reff{linbt}). The condition
that
  \[ \bt\tty(K_j \cd (e_i \ot e_k)) = \bt\tty(e_i \ot e_k)
                       \quad\mbox{for all $j \in J$ and all $i,k \in I$} \]
is satisfied if and only if
  \be \bt\tty(e_i \ot e_k) = 0
    \quad\mbox{for all $i,k \in I$ with $i + k \neq 0$} \;. \label{Kinv} \ee
In particular, this implies that $b$ must be homogeneous of degree zero.

Taking Eqn.~\reff{Kinv} for granted, the conditions
  \[ \bt\tty(E_j \cd (e_i \ot e_k)) = 0 
                       \quad\mbox{for all $j \in J$ and all $i,k \in I$} \]
and
  \[ \bt\tty(F_j \cd (e_i \ot e_k)) = 0 
                       \quad\mbox{for all $j \in J$ and all $i,k \in I$} \]
both yield the same system of linear equations for the elements
$b\tty(e_i\tty, e_k)$. This system has a unique (up to scalar multiples)
solution. Choosing a suitable normalization, the invariant bilinear form
$b^q$ we are looking for is given by
 \[ b^{q}(e_i\tty, e_k) = C^q_{ik} \quad\mbox{for all $i,k \in I$} \,, \]
where
 \[ C^q_{i,\tty k} = C^q_{i,-i} \, \de_{i,-k}
                                     \quad\mbox{for all $i,k \in I$} \,, \]
and where the coefficients $C^q_{i,-i}$ are given by
  \[ C^q_{i,-i} = \left\{
       \begin{array}{cl}
        - q^i & \mbox{for $-1 \geq i \geq -n$} \\[1.0ex]
       q^{-i-2n-2} & \mbox{for $-n-1 \geq i \geq -r$} \\[1.0ex]
          q^i & \mbox{for $1 \leq i \leq n$} \\[1.0ex]
       q^{2n-i} & \mbox{for $n+1 \leq i \leq r$} \;.
       \end{array} \right. \]
Obviously, the matrix $C^q = (C^q_{ij})_{i,\ty j \in I}$ is invertible, i.e.,
the bilinear form $b^q$ is non--degenerate. We note that $C^{q = 1} = C$
(see Eqn.~\reff{defC}), moreover, we have
 \be C^q_{i,-i} \tty C^q_{-i,\ty i} = \left\{
       \begin{array}{rl}
        -1 & \mbox{for $i \in I_{\ol{0}}$} \\[1.0ex]
       q^{-2} & \mbox{for $i \in I_{\ol{1}}$} \;. 
       \end{array} \right. \label{Cqsq} \ee
Thus the matrix $(C^q)^2$ is not equal to $-G$ (recall the
Eqns.~\reff{defG}, \reff{Csq}).

\sect{The structure of the module $V \ot V$ \vspace{-1ex}}
We now are ready to tackle a crucial intermediate problem, namely, to
determine the structure of the tensorial square of the vector module $V$ of
$\Uq$. In the undeformed case, this structure is known. It turns out that
in the deformed case, the structure is completely analogous. In particular,
for $n = m\ty\ty$, the module $V \ot V$ is not completely reducible.
(Actually, if adequately interpreted, the investigations of the present
section apply also to the case $q = 1\ty$.)

To begin with, we stress that the $\Uq$--module $V \ot V$ has the same
weights (with the same multiplicities) as in the undeformed case: For all
$i,j \in I$, the tensor $e_i \ot e_j$ has the weight $\ve_i + \ve_j\ty\ty$.

As expected, $V \ot V$ contains a unique (up to scalar multiples) 
$\Uq$--invariant element, i.e., a non--zero element $a$ such that
 \[ X \cd a = \ve(X) \ty\ty a \quad\mbox{for all $X \in \Uq$} \;. \]
The invariance of $a$ under the action of the generators $K_j$\,; $j \in J$,
is equivalent to the fact that $a$ has the weight zero, i.e., that $a$ is a
linear combination of the following form
  \[ a = \sum_{i \in I} c_i \, e_i \ot e_{-i} \;, \]
with some coefficients $c_i$\ty\ty, $i \in I$.

For an element $a$ of this form, the conditions
 \[ E_j \cd a = 0 \quad\mbox{for all $j \in J$} \vspace{-0.5ex} \]
and
 \[ F_j \cd a = 0 \quad\mbox{for all $j \in J$} \vspace{0.5ex} \]
both yield the same system of linear equations for the coefficients $c_i$\,.
This system has a unique (up to scalar multiples) solution. Choosing a
suitable normalization, the element $a$ is given by
  \[ a = \sum_{i,k \in I} ((C^q)^{-1})_{ik} \, e_i \ot e_k \;, \]
where $C^q$ is the matrix found in Section 4. Of course, this result might
have been anticipated. More explicitly, we have
  \[ a \,=\, \sum_{i=1}^n \tty (q^{-i} \ty e_{-i} \ot e_i
                           - q^{i} \ty e_i \ot e_{-i})
            \,+\, q \!\!\sum_{i=n+1}^r \!(q^{i-2n-1} \ty e_{-i} \ot e_i
                                      + q^{-i+2n+1} \ty e_i \ot e_{-i}) \;. \] 

It is useful to calculate $\bt^q(a)$, where $\bt^q$ is the linear form on
$V \ot V$ defined in Section 4. Setting
  \be d = n - m \;, \label{defd} \ee
we obtain
  \be \bt^q(a)
    \,=\, \frac{1}{q^2 - 1}((q^{-2d} -1) - q^2(q^{2d} - 1))
    \,=\, -\frac{q^d - q^{-d}}{q - q^{-1}}(q^{d+1} + q^{-d-1}) \;.
                                                           \label{bofa} \ee
Note that $\bt^q(a) = 0$ if $n = m\ty\ty$. This is a first indication that
there will be problems in the case $n = m\ty\ty$.

The rest of the present section will now be devoted to prove the following
statements. \\[1.0ex]
a) The $\Uq$--module $V \ot V$ is the direct sum of two submodules
$(V \ot V)_s$ and $(V \ot V)_a\:$:
 \[ V \ot V = (V \ot V)_s \op (V \ot V)_a \;, \]
where in the undeformed case $(V \ot V)_s$ corresponds to the subspace of
super--symmetric and $(V \ot V)_a$ to the subspace of super--skew--symmetric
tensors in $V \ot V$. \\[1.0ex]
b) The submodule $(V \ot V)_s$ is irreducible. \\[1.0ex]
c) The submodule $(V \ot V)_a$ contains a submodule $(V \ot V)_a^0$ of
codimension one, and $(V \ot V)_s \op (V \ot V)_a^0$ is the kernel of the
linear form $\bt^q$ found in Section 4. \\[1.0ex]
d) If $n \neq m\tty$, the $\Uq$--module $(V \ot V)_a^0$ is irreducible, and
$(V \ot V)_a$ is the direct sum of the submodules $(V \ot V)_a^0$ and
$\CC\,a$\,:
 \[ (V \ot V)_a = (V \ot V)_a^0 \op \CC\,a \quad\mbox{if $n \neq m$} \;. \]
e) If $n = m\tty$, we have
 \[ a \in (V \ot V)_a^0 \quad\mbox{if $n = m$} \;, \]
and $(V \ot V)_a^0$ does not have a module complement in $(V \ot V)_a\ty\ty$.
                                                                    \\[1.0ex]
f) If $n = m \geq 2\ty\ty$, the $\Uq$--module $(V \ot V)_a^0/\CC\,a$ is
irreducible. \\[1.0ex]
g) For $n = m = 1\ty$, there exist two submodules $V_4$ and $\Vf$ of
$(V \ot V)_a^0$ such that
 \[ V_4 + \Vf = (V \ot V)_a^0  \quad,\quad V_4 \cap \Vf = \CC\,a \;,\]
and the modules $V_4/\CC\,a$ and $\Vf/\CC\,a$ are irreducible. \\[1.0ex]
h) Let $P_s$ be the projector of $V \ot V$ onto $(V \ot V)_s$ with kernel
$(V \ot V)_a\ty\ty$, and let
 \[ K : V \ot V \lra V \ot V \]
be the linear map defined by
 \be K(u) = \bt^q(u)\ty\ty a \quad\mbox{for all $u \in V \ot V$} \,.
                                                            \label{defK} \ee
Then $\id_{V \ot V}\ty\ty$, $P_s\ty\ty$, and $K$ form a basis of the space
of all $\Uq$--module endomorphisms of $V \ot V$.

In the proof of these claims, we shall obtain more detailed information on
the submodules mentioned above. In particular, we shall construct bases of
the vector spaces $(V \ot V)_s\ty\ty$, $(V \ot V)_a^0\ty\ty$, and
$(V \ot V)_a\ty\ty$.

\subsection{The module $(V \ot V)_s$ \vspace{-1ex}}
As already mentioned above, in the undeformed case the module $(V \ot V)_s$
corresponds to the subspace of all super--symmetric tensors in $V \ot V$.
This subspace is an irreducible $\spob$--submodule of $V \ot V$ and has the
highest weight $\ve_{-r} + \ve_{-r+1}\tty$. If there exists a corresponding
primitive vector in the $\Uq$--module $V \ot V$, it must be a linear 
combination of $e_{-r} \ot e_{-r+1}$ and $e_{-r+1} \ot e_{-r}$\,. Indeed,
there is a \mbox{unique} (up to scalar multiples) linear combination of these
tensors that is annihilated by all $E_j\ty\ty$, $j \in J$. Choosing a
suitable normalization, it is equal to
  \[ s_{-r,-r+1} = e_{-r} \ot e_{-r+1}
                 + \si_{-r,-r+1} \tty q^{-1} \ty e_{-r+1} \ot e_{-r} \;. \]
By definition, $(V \ot V)_s$ is the submodule of $V \ot V$ generated by
$s_{-r,-r+1}\tty$.

Let us define the following elements of $V \ot V$\,:
  \be \begin{array}{rcll}
    s_{i,\ty j} \eqq e_i \ot e_j + \si_{i,\ty j} \ty\ty q^{-1} \ty e_j \ot e_i
              & \mbox{for $i,j \in I$; $i < j$ but $i \neq -j$}  \\[1.0ex]
    s_{i,\ty i} \eqq e_i \ot e_i & \mbox{for $i \in I_{\bar{0}}$} \;,
     \end{array} \label{sij} \ee
furthermore,
  \be s_1 \,=\, e_{-1} \ot e_1 + q^{-2} \ty e_1 \ot e_{-1} \;, \label{so}\ee
and for $2 \leq j \leq r$
  \be \begin{array}{rcl}
     s_j \eqq q^{-\si_{j-1}} \ty e_{-j+1} \ot e_{j-1}
              + \si_{j-1} \tty q^{-1} \ty e_{j-1} \ot e_{-j+1} \\[0.8ex]
         &  & \!\!-\, \si_{j-1} \ty\ty \si_j \ty\ty e_{-j} \ot e_j
              -\si_{j-1} \tty q^{-1} \tty q^{-\si_j} \ty e_j \ot e_{-j} \;.
     \end{array} \label{str} \ee
It turns out that the tensors \reff{sij}, \reff{so}, and \reff{str} form a
basis of $(V \ot V)_s\ty\ty$.

First of all, one shows that the tensors \reff{sij}, \reff{so}, and \reff{str}
can be obtained by iterated action of the $F$--generators on
$s_{-r,-r+1}\tty$. Next one proves that the vector space spanned by these
tensors is a $\Uq$--submodule of $V \ot V$. Since these tensors obviously
are linearly independent, this implies our claim.

Next we have to show that the $\Uq$--module $(V \ot V)_s$ is irreducible.
This can be done as follows. First, we prove the following statement:
                                                                  \\[1.0ex]
If $x$ is a non--zero element of $(V \ot V)_s\ty\ty$, there exists a
monomial $P$ in the $F$--generators such that $P \cd x$ is a non--zero
scalar multiple of $s_{r-1,\ty r}\ty\ty$.

A moment's thought shows that this is a consequence of the following fact:
                                                                  \\[1.0ex]
If $x$ is a (non--zero) weight vector of $(V \ot V)_s$ whose weight is
different from $\ve_{r-1} + \ve_r$ (i.e., if $x$ is not a scalar multiple
of $s_{r-1,\ty r}$\ty), there exists an index $j \in J$ such that
$F_j \cd x \neq 0$\,.

Since $s_{-r,-r+1}\tty$ is a cyclic vector of the $\Uq$--module
$(V \ot V)_s\ty\ty$, the irreducibility of this module will follow if we
can show that there exists a monomial in the $E$--generators that maps
$s_{r-1,\ty r}$ onto a non--zero scalar multiple of $s_{-r,-r+1}\tty$. 
Similar as above, this is a consequence of the following fact:
                                                                  \\[1.0ex]
If $x$ is a (non--zero) weight vector of $(V \ot V)_s$ whose weight is
different from $\ve_{-r} + \ve_{-r+1}$ (i.e., if $x$ is not a scalar multiple
$s_{-r,-r+1}$\ty), there exists an index $j \in J$ such that
$E_j \cd x \neq 0$\,.

The proof of the foregoing statements amounts to easy but lengthy
calculations. Let us mention that one needs the intermediate result that
 \[ a \notin (V \ot V)_s \;. \]
Summarizing, we have proved that the statement b) above is correct.

Obviously, we have
 \[ \bt^q(s_{-r,-r+1}) = 0 \;. \]
Since the $\Uq$--module $(V \ot V)_s$ is irreducible, we conclude that it
is contained in the kernel of $\bt^q$. This proves part of statement c).

\subsection{The module $(V \ot V)_a^0$ \vspace{-1ex}}
Basically, our treatment of the $\Uq$--module $(V \ot V)_a^0$ follows
similar lines to that of $(V \ot V)_s\ty\ty$, however, in the cases $n = m$
there are several complications.

In the undeformed case, the module $(V \ot V)_a^0$ corresponds to the
subspace of all super--skew--symmetric tensors in $V \ot V$ with a
vanishing ``symplectic trace'' (i.e., which belong to the kernel of
$\bt\tty$). As an $\spob$--submodule of $V \ot V$, it is generated by the
tensors $e_{-r} \ot e_{-r}$ and $e_r \ot e_r\ty\ty$, and for
$(n\ty,m) \neq (1\ty,1)$, each of these tensors alone is already sufficient.

In the present deformed case, it is easy to see that $e_{-r} \ot e_{-r}$ is
annihilated by the $E$--generators, and that $e_r \ot e_r\ty\ty$ is
annihilated by the $F$--generators. Accordingly, we define $(V \ot V)_a^0$
to be the $\Uq$--submodule of $V \ot V$ generated by $e_{-r} \ot e_{-r}$
and $e_r \ot e_r\ty\ty$.

Let us define the following elements of $V \ot V$\,:
  \be \begin{array}{rcll}
    a_{i,\ty j} \eqq e_i \ot e_j - \si_{i,\ty j} \tty\ty q \tty\ty e_j \ot e_i
              & \mbox{for $i,j \in I$; $i < j$ but $i \neq -j$}  \\[1.0ex]
    a_{i,\ty i} \eqq e_i \ot e_i & \mbox{for $i \in I_{\bar{1}}$} \;,
     \end{array} \label{aij} \ee
and for $2 \leq j \leq r$
  \be \begin{array}{rcl}
     a_j \eqq q^{-\si_{j-1}} \ty e_{-j+1} \ot e_{j-1}
              - \si_{j-1} \ty\ty q \tty\ty e_{j-1} \ot e_{-j+1} \\[0.8ex]
         &  & \!\!-\, \si_{j-1} \ty\ty \si_j \ty\ty e_{-j} \ot e_j
              +\si_{j-1} \ty\ty q \tty\ty q^{-\si_j} \ty e_j \ot e_{-j} \;.
     \end{array} \label{atr} \ee
(I hope there is no risk to confound the tensors $a_{i,\ty j}$ with the
\vspace{0.2ex}
elements of the Cartan matrix.) It turns out that the tensors \reff{aij} and
\reff{atr} form a basis of $(V \ot V)_a^0\ty\ty$.

First one proves that the vector space $U$ spanned by the tensors \reff{aij}
and \reff{atr} is a $\Uq$--submodule of $V \ot V$. Since these tensors
obviously are linearly independent, they form a basis of $U$.

Next one shows that $a_{-r,-r}$ generates this module, {\em provided} that
$r \geq 3$ (i.e., provided that $(n\ty,m) \neq (1\ty,1)$\ty). In the case
$n = m = 1\ty$, we have
 \[ a_2 = a \quad\mbox{if $m = n = 1$} \;, \]
and $a_{-2,-2}$ generates the $U_q(\mathfrak{spo}(2 \ty|\ty 2))$--submodule
 \[ V_4 \,=\, \CC\,a_{-2,-2}
                 \op \CC\,a_{-2,-1} \op \CC\,a_{-2,\ty 1} \op \CC\,a \;, \]
while $a_{2,\ty 2}$ generates the $U_q(\mathfrak{spo}(2 \ty|\ty 2))$--submodule
 \[ \Vf \,=\, \CC\,a_{2,\ty 2} \op \CC\,a_{1,\ty 2}
                                    \op \CC\,a_{-1,\ty 2} \op \CC\,a \;. \]
Obviously, we have
  \[ V_4 + \Vf = U \quad,\quad V_4 \cap \Vf = \CC\,a \;, \]
and it is easy to see that the $U_q(\mathfrak{spo}(2 \ty|\ty 2))$--modules
$V_4 / \CC\ty\ty a$ and $\Vf / \CC\ty\ty a$ are irreducible. This proves
statement g), and we also have shown that in all cases
 \[ U = (V \ot V)_a^0 \;. \]

Let us next show that the sum of the subspaces $(V \ot V)_s$ and
$(V \ot V)_a^0$ of $V \ot V$ is direct. Obviously, it is sufficient to prove
the analogous statement for the corresponding weight spaces. For non--zero
weights, this is trivial. To prove the claim for the weight zero, it is
sufficient to show that the tensors $s_i$\tty, $1 \leq i \leq r\ty$, and
$a_j$\tty, $2 \leq j \leq r\ty$, are linearly independent. To show this, we
observe that
  \be s_j = u_j + q^{-1} \ty v_j \quad,\quad
     a_j = u_j - q \tty v_j \;\;\quad\mbox{for $2 \leq j \leq r$} \;,
                                                          \label{intuv} \ee
where the tensors $u_j$ and $v_j\ty\ty$; $2 \leq j \leq r\ty$, are defined by
  \[ \begin{array}{rcl}
      u_j \eqq q^{-\si_{j-1}} \ty e_{-j+1} \ot e_{j-1} 
                   - \si_{j-1} \ty\ty \si_j \ty\ty e_{-j} \ot e_j \\[0.8ex]
      v_j \eqq  \si_{j-1} \ty\ty e_{j-1} \ot e_{-j+1}
                   - \si_{j-1} \ty\ty q^{-\si_j} \ty e_j \ot e_{-j} \;.
     \end{array} \]
Consequently, we have to prove that the tensors $s_1$ and $u_j$\tty,
$v_j$\tty, $2 \leq j \leq r\ty$, are linearly independent. This follows
from the obvious fact that the $2r$ tensors $e_{-1} \ot e_1\ty$, $s_1$\ty,
\mbox{$u_j$\tty,} $v_j$ span the same subspace of $V \ot V$ as the $2r$
tensors $e_i \ot e_{-i}\ty\ty$, $i \in I$, namely, the weight space of
$V \ot V$ corresponding to the weight zero.

The proof above shows that the codimension of $(V \ot V)_s \op (V \ot V)_a^0$
in $V \ot V$ is equal to one. Obviously, $\bt^q$ vanishes on the tensors
\vspace{0.2ex}
$e_{-r} \ot e_{-r}$ and $e_r \ot e_r\ty\ty$, hence also on the submodule
$(V \ot V)_a^0$ generated by them. As noted earlier, $\bt^q$ also vanishes
on $(V \ot V)_s\ty\ty$. Since $\bt^q$ is a non--zero linear form on
$V \ot V$, it follows that $(V \ot V)_s \op (V \ot V)_a^0$ is the kernel of
$\bt^q$. This proves the last claim of statement c).

Using Eqn.~\reff{bofa}, our last result implies that
 \be a \in (V \ot V)_s \op (V \ot V)_a^0
                 \quad\mbox{if and only if $n = m$} \;. \label{ainkerbt} \ee
Since the $\Uq$--module $(V \ot V)_s$ is irreducible and not
one--dimensional, it follows that
 \be a \in (V \ot V)_a^0 \quad\mbox{for $n = m$} \;. \label{ainvvaz} \ee
Indeed, it can be shown that
  \be a \,=\, \sum_{j=2}^{n+1}\tty [\tty j - 1] \tty a_j
                       -\sum_{j=n+2}^{2n} [\ty 2n + 1 - j\ty] \tty a_j
                               \quad\mbox{for $n = m$} \;, \label{aneqm} \ee
where, for all integers $s\tty$, the $q$--number $[s]$ is defined by
  \[ [s] \,=\, [s]_q \,=\, \frac{q^s - q^{-s}}{q - q^{-1}} \;. \]
For $n = 1\ty$, Eqn.~\reff{aneqm} is just the equation $a = a_2$ mentioned
earlier.

Using Eqn.~\reff{ainkerbt} and the fact that $(V \ot V)_s \op (V \ot V)_a^0$
is equal to the kernel of $\bt^q$, it follows that, for $n = m\ty\ty$, this
submodule does not have a module complement in $V \ot V$. In fact, any such
complement would have to be a trivial one--dimensional submodule of
$V \ot V$, and hence would be spanned by an invariant element of $V \ot V$.
But the invariant elements of $V \ot V$ are the scalar multiples of
$a\ty\ty$. Combined with Eqn.~\reff{ainvvaz}, this yields statement e).

Finally, to answer questions of irreducibility, we prove the following
technical results. \\[1.0ex]
Suppose that $r \geq 3\ty\ty$, and that $x \in (V \ot V)_a^0$ is a
(non--zero) weight vector which is neither proportional to $a_{r,\ty r}$ nor
to $a\tty$. Then there exists an index $j \in J$ such that
$F_j \cd x \notin \CC\,a$ (in particular, we have $F_j \cd x \neq 0$\ty).
                                                                   \\[1.0ex]
Suppose that $r \geq 3\ty\ty$, and that $x \in (V \ot V)_a^0$ is a
(non--zero) weight vector which is neither proportional to $a_{-r,-r}$ nor
to $a\tty$. Then there exists an index $j \in J$ such that
$E_j \cd x \notin \CC\,a$ (in particular, we have $E_j \cd x \neq 0$\ty).

As in the case of $(V \ot V)_s\ty\ty$, these results follow from easy but
lengthy calculations. Once they are established, it is easy to draw the
following conclusions. \\[1.0ex]
If $n \neq m\ty\ty$, the $\Uq$--module $(V \ot V)_a^0$ is irreducible.
                                                                   \\[1.0ex]
If $n = m \geq 2\ty\ty$, the $\Uq$--module $(V \ot V)_a^0/\CC\,a$ is
irreducible. \\[1.0ex]
These results prove statement f) and the first claim of statement d).

\subsection{The module $(V \ot V)_a$ \vspace{-1ex}}
Our next task is to construct the submodule $(V \ot V)_a$ of $V \ot V$. In
the case $n \neq m\tty$, this is easy. Recalling Eqn.~\reff{ainkerbt} and
the fact that $(V \ot V)_s \op (V \ot V)_a^0$ has the codimension one in
$V \ot V$, it follows that
  \be V \ot V \,=\, (V \ot V)_s \op (V \ot V)^0_a \op \CC\,a
                               \quad\mbox{if $n \neq m$} \;. \label{dec} \ee
As we know, the submodules on the right hand side of this equation are
irreducible, moreover, they are obviously non--isomorphic. This implies that
Eqn.~\reff{dec} is the unique decomposition of the $\Uq$--module $V \ot V$
into irreducible submodules. Setting
  \[ (V \ot V)_a = (V \ot V)_a^0 \op \CC\,a \;, \]
we have proved the statements a)--d) in the case $n \neq m\tty$.

Unfortunately, this type of reasoning is not possible in the case
$n = m\tty$. Since we want to obtain a unified treatment of the problem, we
start all over again and present an approach which is applicable in all
cases.

To begin with, we note that every module complement of $(V \ot V)_s$ in
$V \ot V$ must take the form $(V \ot V)_a^0 \op \CC\,g\ty\ty$, where $g$ is
a weight vector of $V \ot V$ of weight zero. Indeed, since the tensors
$e_{-r} \ot e_{-r}$ and $e_r \ot e_r$ do not belong to $(V \ot V)_s\ty\ty$,
and since the corresponding weights have multiplicity one, these two tensors
and hence the submodule generated by them must be contained in every module
complement. Obviously, a subspace of this type is a module complement of
$(V \ot V)_s$ if and only if $g \notin (V \ot V)_s \op (V \ot V)_a^0$ and if
$E_j \cd g$ and $F_j \cd g$ belong to $(V \ot V)_a^0\ty\ty$, for all
$j \in J$.

Since $g$ is of zero weight, it takes the form
  \be g \,=\, \sum_{i \in I} g_i \, e_i \ot e_{-i} \;, \label{ggen} \ee
with some coefficients $g_i \in \CC$\,. If $E_j \cd g$ is non--zero, it is
a weight vector with non--zero weight and hence belongs to
$(V \ot V)_s \op (V \ot V)_a^0\ty\ty$. Consequently, $E_j \cd g$ lies in
$(V \ot V)_a^0$ if and only if its component in $(V \ot V)_s$ is equal to
zero. It follows that we have $E_j \cd g \in (V \ot V)_a^0$ if and only if
the coefficients $g_i$ satisfy the following system of linear equations:
  \be g_1 + q^2 \ty g_{-1} = 0 \label{go} \ee
  \be q \tty g_{j+1} - q \tty q^{\si_j}\ty g_j
     \,=\, \si_j \tty g_{-j} - \si_{j+1} \tty q^{\si_{j+1}}\ty g_{-j-1}
                   \quad\mbox{for $-r \leq j \leq -2$} \;. \label{gtr} \ee
The condition that $F_j \cd g \in (V \ot V)_a^0$ for all $j \in J$ is
equivalent to the same system of equations.

The general solution of this system can easily be described: We can choose
$g_{-1}$\ty, $g_{-2}$\tty,\tty\tty\ldots, $g_{-r}$ arbitrarily, and then the
coefficients $g_1$\ty, $g_2$\tty,\tty\tty\ldots, $g_r$ are uniquely fixed.

Let $X_a$ be the subspace of $V \ot V$ consisting of all tensors of the form
\reff{ggen} such that the coefficients $g_i$ satisfy the system \reff{go},
\reff{gtr}. According to the foregoing result, this subspace is
$r$--dimensional. Obviously, $X_a$ contains the $(r - 1)$--dimensional
weight space of $(V \ot V)_a^0$ corresponding to the weight zero. On the
other hand, $X_a$ does not contain any non--zero elements of
$(V \ot V)_s\ty\ty$ (indeed, any such element would be invariant and hence
proportional to $a$\ty). It follows that
 \[ (V \ot V)_a = (V \ot V)_a^0 + X_a \]
is a module complement of $(V \ot V)_s$ in $V \ot V$.

By this latter property, $(V \ot V)_a$ is uniquely fixed. In fact, let
$(V \ot V)'_a$ be an arbitrary module complement of $(V \ot V)_s$ in
$V \ot V$. As noted at the beginning of this discussion, $(V \ot V)'_a$
contains the submodule $(V \ot V)_a^0\ty\ty$. Let $X'_a$ be the weight space
of $(V \ot V)'_a$ corresponding to the weight zero. Of course, $X'_a$ is
$r$--dimensional. Consider an arbitrary element $g \in X'_a$. For every
$j \in J$, the elements $E_j \cd g$ and $F_j \cd g$ lie in $(V \ot V)'_a$
and have a non--zero weight, which implies that they are elements of
$(V \ot V)_a^0\ty\ty$. By the definition of $X_a\ty\ty$, this proves that
$g \in X_a\ty\ty$. Thus we have shown that $X'_a \subset X_a\ty\ty$, and for
reasons of dimension, this implies that $X'_a = X_a\ty\ty$. It follows that
$X_a \subset (V \ot V)'_a$\,, hence that
  \[ (V \ot V)_a \subset (V \ot V)'_a \;, \]
and finally 
  \[ (V \ot V)_a = (V \ot V)'_a \;, \]
as claimed.

We proceed by choosing, in a unified way, an element $t \in X_a$ that does
not belong to $(V \ot V)_a^0\ty\ty$. Let $t$ be the tensor of the form
\reff{ggen} whose coefficients $g^{\tty t}_i$ are fixed by the requirement
that
  \[ g^{\tty t}_{-1} = 1 \quad,\quad g^{\tty t}_j = 0
                                   \;\;\mbox{for $-r \leq j \leq -2$} \;. \]
Recall that the coefficients $g^{\tty t}_i$ with $1 \leq i \leq r$ can then
be calculated by means of the system \reff{go}, \reff{gtr}. We obtain
 \be t \,=\, e_{-1} \ot e_1 - q^2 \tty e_1 \ot e_{-1}
                 + (q - q^{-1})\! \sum_{i=2}^{n+m} \ty
                      \Cqi_{i,-i} \tty\ty e_i \ot e_{-i} \;. \label{deft} \ee
It is easy to check that
  \[ \bt^q(t) = -q^{-1} \ty (q^{2d+2} + 1) \]
(where $d = n - m\ty\ty$; see Eqn.~\reff{defd}). Thus $t$ does not belong to
the kernel of $\bt^q$ (since $q$ is not a root of unity).

According to the results of the present subsection, the tensor $t\tty$ and
the tensors \reff{aij} and \reff{atr} form a basis of the vector space
$(V \ot V)_a\ty\ty$. \vspace{2.0ex}

\noindent
{\em Remark 5.1.} We note that some other, simpler looking choices for the
tensor $t$ are possible. For example, the tensor
$e_{-r} \ot e_r + e_r \ot e_{-r}$ is a candidate. However, it is not at all
obvious that such a choice would simplify the subsequent
\vspace{2.0ex}
calculations.

The reader will easily convice himself/herself that, at this stage, we have
proved the statements a)--g).

\subsection{The module endomorphisms of $V \ot V$ \vspace{-1ex}}
In the present subsection we are going to prove statement h), i.e., that the
linear operators $\id_{V \ot V}\ty\ty$, $P_s\ty\ty$, and $K$ form a basis of
the space of all endomorphisms of the $\Uq$--module $V \ot V$ (recall that
$P_s$ is the projector of $V \ot V$ onto $(V \ot V)_s$ with kernel
$(V \ot V)_a\ty\ty$, and that the map $K$ has been defined in
Eqn.~\reff{defK}). Obviously, the maps $\id_{V \ot V}$ and $P_s\ty\ty$ are
module endomorphisms, and since $\bt^q$ and $a$ are invariant, the same is
true of $K\ty\ty$.

Now let $Q$ be a module endomorphism of $V \ot V$, i.e., an even linear map
of $V \ot V$ into itself that commutes with the action of $\Uq$. Since
every weight vector of $V \ot V$ of weight $\ve_{-r} + \ve_{-r+1}$, that is
annihilated by all $E_j\ty\ty$, $j \in J$, is proportional to
$s_{-r,-r+1\tty}$, there exists a constant $c_s$ such that
  \[ Q(s_{-r,-r+1}) = c_s\tty s_{-r,-r+1} \;. \]
Similarly, since the multiplicity of the weight $2\tty\ve_{-r}$ is equal to
one, we have
  \[ Q(a_{-r,-r}) = c_a\tty a_{-r,-r} \;, \]
with some constant $c_a\ty\ty$. It follows that
  \be \begin{array}{rcll}
     Q(x) \eq c_s\tty x & \;\;\mbox{for all $x \in (V \ot V)_s$} \\[0.7ex]
     Q(y) \eq c_a\tty y & \;\;\mbox{for all $y \in (V \ot V)^0_a$} \;.
     \end{array} \label{Qev} \ee
Indeed, since the tensor $s_{-r,-r+1}$ generates the submodule
$(V \ot V)_s$\,, the first of these equations follows immediately.
Similarly, for $r \geq 3$ the tensor $a_{-r,-r}$ generates the submodule
$(V \ot V)^0_a$\,, which implies the second equation in this case.

In the case $r = 2$\,, i.e., for $m = n = 1$\,, we argue as follows. Quite
generally, we also have
  \[ Q(a_{r,\ty r}) = \ol{c}_a\tty a_{r,\ty r} \;, \]
with some constant $\ol{c}_a\ty\ty$. Then, for $m = n = 1\ty$, we conclude
as above that
  \[ \begin{array}{rcll}
        Q(y) \eq c_a\tty y & \;\;\mbox{for all $y \in V_4$} \\[0.7ex]
   Q(\ol{y}) \eq \ol{c}_a\tty \ol{y}
                         & \;\;\mbox{for all $\ol{y} \in \Vf$} \;.
     \end{array} \]
Since $a$ lies in $V_4$ and in $\Vf\ty\ty$, this implies that
  \[ c_a = \ol{c}_a \;, \]
and since
  \[ V_4 + \Vf = (V \ot V)^0_a \;, \]
it follows that the second of the Eqns.~\reff{Qev} holds for $n = m = 1$ as
well.

Now let
  \[ P_a =  \id_{V \ot V} - P_s \]
be the projector of $V \ot V$ onto $(V \ot V)_a$ with kernel
$(V \ot V)_s\ty\ty$. Then the equations derived above show that
  \[ Q' = Q - c_s\tty P_s - c_a\tty P_a \]
is a module endomorphism of $V \ot V$, that vanishes on
  \[ \mbox{ker}\,\bt^q = (V \ot V)_s \op (V \ot V)^0_a \;. \]
Consequently, it induces a module homomorphism
  \[ (V \ot V)/\ty\mbox{ker}\,\bt^q \lra V \ot V \,. \]
The module on the left hand side is one--dimensional and trivial, and all
invariants in $V \ot V$ are proportional to $a\ty\ty$. This implies that
  \[ Q' = c_0\tty K \;, \]
with some constant $c_0\ty\ty$, which proves our claim.

We close this subsection by the remark that the maps $\id_{V \ot V}\ty\ty$,
$P_s\ty\ty$, and $K$ commute one with another.

\sect{Calculation of the $R$--matrix \vspace{-1ex}}
At last, we are prepared to calculate the $R$--matrix $R$ or, equivalently,
the braid generator $\Rh$ of $\Uq$ in the vector representation. By
definition, $R$ is the representative of the universal $R$--matrix $\cal R$
in the vector representation,
 \[ R = {\cal R}_{V \ot V} \;, \]
and $\Rh$ is given by
 \[ \Rh = P R \;, \]
where
 \[ P : V \ot V \lra V \ot V \]
denotes the twist operator (in the super sense), which is given by
 \[ P(x \ot y) = (-1)^{\xi\eta} \ty y \ot x \;, \]
for all $x \in V_{\xi}\ty\ty$, $y \in V_{\eta}\ty\ty$, with
$\xi,\eta \in \ZZ_2\,$. To calculate the $R$--matrix (or the braid generator)
means to calculate its matrix elements with respect to the basis \\
$(e_i \ot e_j)_{i,j \in I}$ of $V \ot V$. \vspace{2.0ex}

\noindent
{\em Remark 6.1.} Due to the fact that $\cal R$ is given in terms of a
formal power series, the foregoing remarks need to be amended. See below
for further details. \vspace{2.0ex}

In order to perform the calculation we observe that $\Rh$ is an endomorphism
of the $\Uq$--module $V \ot V$. According to Section 5.4, this implies that
$\Rh$ is a linear combination of $\id_{V \ot V}\ty\ty$, $K\ty$, and
$P_s\ty\ty$. Since the matrix elements of $\id_{V \ot V}$ and $K\ty$ are
known, our task consists of two pieces, namely, to calculate the projector
$P_s$ (or any module endomorphism of $V \ot V$ ``containing'' $P_s$ in a
non--trivial way), and to find the aforementioned linear combination. As we
are going to see, the second problem can easily be dealt with once the first
problem has been solved.

Obviously, every module endomorphism of $V \ot V$ maps each of the weight
spaces into itself. Consequently, the first problem splits into a number of
subproblems, one for each of the weight spaces of $V \ot V$. Since the weight
spaces corresponding to the non--zero weights are at most two--dimensional,
the corresponding subproblems are trivial, and we are left with the
subproblem corresponding to the zero weight. Basically, this latter
problem amounts to writing the tensors $e_i \ot e_{-i}\,$, $i \in I$, as
linear combinations of the tensors \reff{sij}\,--\,\reff{atr} and $t\ty\ty$,
i.e., we have to invert a certain $(2\ty r \ti 2\ty r)$--matrix, whose
elements are rational functions of $q\ty\ty$.

Unfortunately, the corresponding calculations turn out to be rather tedious.
Accordingly, I don't present the details of this calculation but only mention
two points. First, in the course of the calculations I have taken advantage
of the tensors $u_j$ and $v_j$ introduced in Eqn.~\reff{intuv} and of the
resulting equations
  \[ P_s(u_j) \,=\, q \tty (q + q^{-1})^{-1} \ty s_j \quad,\quad
                            P_s(v_j) \,=\, (q + q^{-1})^{-1} \ty s_j \;. \]
Secondly, I have applied the following simple trick. In Ref.~\cite{RTF},
the formulae for $\Rh$ are slightly simpler than those for $P_s\ty\ty$.  On
the other hand, again according to Ref.~\cite{RTF} (see also
Ref.~\cite{ZBG}), it is tempting to conjecture that
 \be \Rh = \Rh\ty' \;, \label{Rhcon} \ee
where the endomorphism $\Rh\ty'$ of $V \ot V$ is defined by
 \be \Rh\ty' = (q + q^{-1})\tty P_s - q^{-1} \II \ot \II
                  - (q - q^{-1})(1 + q^{2d+2})^{-1} K \;. \label{defRp} \ee
Here and in the following, $\II$ denotes the unit operator of $V\ty$:
 \[ \II = \id_V \;. \]
Accordingly, I haven't calculated $P_s$ but rather the operator $\Rh\ty'$.
The Eqn.~\reff{Rhcon} can then be proved at a later stage (which solves the
second problem mentioned at the beginning of this section).

A long calculation shows that for $1 \leq j \leq n$
  \[ \begin{array}{rcl}
      \Rh\ty'(e_{-j} \ot e_j)
        \eq (q - q^{-1})\ty\ty\dis{\sum_{i=1}^n}\,\ty
                                     q^{-j-i} \ty e_{-i} \ot e_i \\[2.1ex]
    \!& &\hspace{-1.2em}
       {} + (q - q^{-1})\!\dis{\sum_{i=n+1}^{n+m}}\!
                                  q^{-2n+i-j} \ty e_{-i} \ot e_i \\[2.1ex]
    \!& &\hspace{-1.2em}
       {} - (q - q^{-1})\ty\ty\dis{\sum_{i=1}^{j-1}}\,
                                      q^{i-j} \ty e_i \ot e_{-i} \\[2.3ex]
    \!& &\hspace{-1.2em}
       {} + (q - q^{-1})\, e_{-j} \ot e_j \tty +\tty q^{-1} \ty e_j \ot e_{-j}
     \end{array} \]
  \[ \begin{array}{rcl}
      \Rh\ty'(e_j \ot e_{-j})
        \eq -(q - q^{-1})\!\dis{\sum_{i=j+1}^n}\,\!
                                  q^{\ty j-i} \ty e_{-i} \ot e_i \\[2.1ex]
    \!& &\hspace{-1.2em}
        {} - (q - q^{-1})\!\dis{\sum_{i=n+1}^{n+m}}\!
                                  q^{-2n+i+j} \ty e_{-i} \ot e_i
                                                  \hspace{1.0em} \\[2.3ex]
    \!& &\hspace{-1.2em}
        {} + q^{-1} \ty e_{-j} \ot e_j \;,
     \end{array} \]
and for $n+1 \leq j \leq n+m$ \vspace{-0.2ex}
  \[ \begin{array}{rcl}
      \Rh\ty'(e_{-j} \ot e_j)
        \eq\! -(q - q^{-1})\ty\ty\dis{\sum_{i=1}^n}\,\ty
                               q^{-2n-i+j-2} \tty e_{-i} \ot e_i \\[2.1ex]
    \!& &\hspace{-1.2em}
        {} - (q - q^{-1})\!\dis{\sum_{i=n+1}^{n+m}}\!
                               q^{-4n+i+j-2} \tty e_{-i} \ot e_i \\[2.1ex]
    \!& &\hspace{-1.2em}
        {} + (q - q^{-1})\ty\ty\dis{\sum_{i=1}^{n}}\,
                               q^{-2n+i+j-2} \tty e_i \ot e_{-i} \\[2.1ex]
    \!& &\hspace{-1.2em}
        {} - (q - q^{-1})\!\dis{\sum_{i=n+1}^{j-1}}\!
                                  q^{\ty j-i} \ty e_i \ot e_{-i} \\[2.5ex]
    \!& &\hspace{-1.2em}
        {} + (q - q^{-1})\, e_{-j} \ot e_j \tty -\tty q \tty e_j \ot e_{-j}
     \end{array} \]
  \[ \begin{array}{rcl}
      \Rh\ty'(e_j \ot e_{-j})
        \eq\! -(q - q^{-1})\!\dis{\sum_{i=j+1}^{n+m}}\!
                                      q^{i-j} \ty e_{-i} \ot e_i
                                                  \hspace{2.5em} \\[2.3ex]
    \!& &\hspace{-1.2em}
        {} - q \tty e_{-j} \ot e_j \;.
     \end{array} \]

It might have been difficult to unify these equations in a concise formula.
Fortunately, Ref.~\cite{RTF} suggests that, for all $i \in I$, we have
  \[ \begin{array}{rcl}
     \Rh\ty'(e_i \ot e_{-i})
       \eq \si_i \, q^{-\si_i} \tty e_{-i} \ot e_i
           \ty +\ty (q - q^{-1})\, \theta(-i > i)\, e_i \ot e_{-i} \\[1.5ex]
    \!& &\hspace{-1.2em}
              \dis{{} - (q - q^{-1}) \, C^q_{i,-i}
                 \sum_{k > i} \Cqi_{-k\ty,\tty k} \, e_{-k} \ot e_k} \;,
     \end{array} \]
where, for all $i,j \in I$, the symbol $\theta(\ty j > i)$ is defined by
  \[ \theta(\ty j > i) = \left\{ \begin{array}{rl}
                             1 & \mbox{if $j > i$} \\
                             0 & \mbox{otherwise\,.}
                             \end{array} \right. \]
It is not difficult to see that this is indeed the case.

The remaining tensors $\Rh\ty'(e_i \ot e_j)\ty$, with $i,j \in I$, are easily
determined: \\[0.5ex]
If $i < j$, but $i \neq -j$, we obtain
  \[ \begin{array}{rcl}
      \Rh\ty'(e_i \ot e_j)
         \eqq \si_{i,\ty j} \tty e_j \ot e_i
                                  + (q - q^{-1}) \tty e_i \ot e_j \\[0.7ex]
      \Rh\ty'(e_j \ot e_i) \eqq \si_{j,i} \tty e_i \ot e_j \;.
     \end{array} \]
On the other hand, we find for all $i \in I$
  \[ \Rh\ty'(e_i \ot e_i) = \si_i \ty\ty q^{\si_i} \ty e_i \ot e_i \;. \] 

Let us next prove Eqn.~\reff{Rhcon}. Needless to say, at this point we have
to make contact with the theory of the universal $R$--matrix. Fortunately,
only very little of that theory is needed. Hence it should be sufficient to
recall a few simple facts. For more details, we refer the reader to
Ref.~\cite{Yam}.

Usually, the theory is formulated in the framework of formal power series in
one indeterminate $h\tty$. In particular, the complex parameter $q$ is
replaced by
 \[ q = e^h \;, \]
and the corresponding quantum superalgebra $U_h$ is a topological Hopf
superalgebra over the ring $\Ch$ of formal power series in $h\tty$. The
Cartan subalgebra $\fh$ of $\spo$ is regarded as a subspace of $U_h\ty$, and
in this sense, the elements $K_i$ are given by
 \[ K_i = \exp(h\ty H_{\al_i}) \quad,\quad H_{\al_i} = d_i \ty H_i \;, \]
where the elements $H_{\lam}$ have been defined in Eqn.~\reff{Hlam}, and the
elements $H_i = H_i^v$ in Eqn.~\reff{Hjv}.

In this setting, instead of $V$ we have to consider the $\Ch$--module $\Vh$
\vspace{0.2ex}
that is obtained from $V$ by an extension of the domain of scalars from
$\CC$ to $\Ch\,$:
 \[ \Vh \,=\, V \otC \Ch  \;. \]
It is well--known that $\Vh$ is a graded $U_h$--module in a natural way, and
it is this module which in the present setting takes the role of the vector
module. More explicitly, the elements $e_i \ot 1 \in \Vh\ty\ty$, $i \in I$,
form a basis of the $\Ch$--module $\Vh\ty\ty$, and it is customary to
identify $e_i \ot 1$ with $e_i\ty\ty$. With this convention, the action of
the generators $E_j\ty\ty$, $F_j\ty\ty$, $K_j\,$; $j \in J$, on the basis
elements is given by the same formulae as in Section 4 (of course, with a
different meaning of $q$\ty).

The tensor product (over $\Ch\tty$) of $\Vh$ with itself is also a graded
$U_h$--module. On the other hand, it is known that
 \[ \Vh \otCh \Vh \simeq (V \otC V) \otC \Ch \]
(as graded $\Ch$--modules), and the elements $e_i \ot e_j \ot 1\ty\ty$;
$i,j \in I$, form a basis of this $\Ch$--module. Once again, we identify
$e_i \ot e_j \ot 1$ with $e_i \ot e_j\ty\ty$, and then the action of the
generators $E_j\ty\ty$, $F_j\ty\ty$, $K_j$ is given by the same formulae as
in Section 5.

Using these conventions, the arguments of Section 5 can be adopted almost
verbatim. Of course, we have to keep in mind that $\Ch$ is not a field but
only a ring. Correspondingly, at various instances we have to observe
that a scalar is not only different from zero but even invertible. Moreover,
the usual concept of an irreducible module is not useful here: If $W$ is a
graded $U_h$--module, then $h W$ is a graded submodule of $W$ and, in
general, different from $W$.

In particular, the $U_h$--module endomorphisms of $\Vh \ot_{\Ch} \Vh$ are
\vspace{0.2ex}
linear combinations of the identity map, $K\ty$, and $\Rh\ty'$ (interpreted
as $\Ch$--linear maps of $(V \ot V) \ot \Ch$ into itself). Now we can prove
Eqn.~\reff{Rhcon} (in the present setting), but let us first complete our
survey. Eqn.~\reff{Rhcon} shows that $\Rh$ depends on $h$ only through $e^h$
(a result which immediately follows by inspection of the formula for the
universal $R$--matrix). In fact, its matrix elements are Laurent polynomials
in $e^h$. Substituting for $e^h$ the complex number $q$ we started with, we
obtain the braid generator that we want to calculate.

Before we can proceed to the calculation proper, we remind the reader of
the general form of the universal $R$--matrix. Using the fact that the
vector spaces $(\fh \ot \fh)^{\ast}$ and $\fh^{\ast} \ot \fh^{\ast}$ are
canonically isomorphic, it is obvious that there exists a unique tensor
$B \in \fh \ot \fh$ such that
  \[ (\lam \ot \mu)(B) = (\lam \ty|\ty \mu)
                 \quad\mbox{for all $\lam\ty , \mu \in \fh^{\ast}$} \;. \]
In terms of this tensor, we have
  \be {\cal R} = e^{h B} \ty (1 \ot 1 + \ldots\,) \;, \label{Runiv} \ee
where the dots stand for an infinite sum of terms of the form $X \ot X'$,
in which $X$ and $X'$ are weight vectors of the quantum superalgebra with
{\em non--zero} (and opposite) weights (see Refs.~\cite{Yam}, \cite{KTo}).

Now we are ready to prove Eqn.~\reff{Rhcon}. According to the preceding
\vspace{0.2ex}
discussion, $\Rh$ is a $\Ch$--linear combination of $\Rh\ty'\ty$,
$\II \ot \II\,$, and $K\ty$ (regarded as $\Ch$--linear maps of
$(V \ot V) \ot \Ch$ into itself). Equivalently, this means that
  \be {\cal R}_{(V \ot V) \ot \Ch} = a \ty P \Rh\ty' + b \ty P + c \ty P K \;,
                                                            \label{Ans} \ee
with some coefficients $a,b,c \in \Ch\ty$. In order to determine these
coefficients, we apply Eqn.~\reff{Ans} to the tensors $e_i \ot e_j$ and keep
only the diagonal terms, i.e., the terms proportional to $e_i \ot e_j\ty\ty$.
On the left hand side, the terms indicated by the dots in Eqn.~\reff{Runiv}
do not contribute, and we are left with
  \[ \begin{array}{rcl}
   (e^{h B})(e_i \ot e_j) \eq q^{(\ve_i \ty|\tty \ve_j)} \ty e_i \ot e_j
                                                                   \\[0.7ex]
           \eq \left\{ \begin{array}{ll}
                e_i \ot e_j & \mbox{if $i \neq j,-j$}            \\[0.5ex]
                q^{\si_i} \tty e_i \ot e_i & \mbox{if $j = i$}   \\[0.5ex]
                q^{-\si_i} \tty e_i \ot e_{-i} & \mbox{if $j = -i$} \;.
                       \end{array} \right.
     \end{array} \]
Using the formulae for $\Rh\ty'(e_i \ot e_j)$ obtained above, the analogous
terms on the right hand side can easily be calculated. Comparing both sides,
we obtain the following equations: \\
For $i \neq j,-j$
  \[ 1 = a \;, \]
for $i = j$
  \[ q^{\si_i} = a \ty\ty q^{\si_i} + b \ty\ty \si_i \;, \]
for $i = -j$
  \[ q^{-\si_i} = a \ty\ty q^{-\si_i} + c \ty\ty \si_i \;. \]
The unique solution of these equations is
  \[ a = 1 \quad,\quad b = c = 0 \;. \]
This implies that
  \[ {\cal R}_{(V \ot V) \ot \Ch} = P \Rh\ty' \;, \]
which proves our claim. \vspace{2.0ex}

\noindent
{\em Remark 6.1.} It is well--known that if $\cal R$ is a universal
$R$--matrix for a Hopf superalgebra $H$, then so is ${\cal R}_{\ty 21}^{-1}$
(we are using the standard notation). Moreover, if $V$ is any graded
$H$--module, and if $\Rh$ is the braid generator in $V \ot V$ with respect to
$\cal R\,$, then $\Rh^{-1}$ is the braid generator in $V \ot V$ with respect
to ${\cal R}_{\ty 21}^{-1}\ty$. Thus we can apply the preceding discussion to
${\cal R}_{\ty 21}^{-1}$ and $\Rh^{-1}$. First of all, we conclude that
 \[ P\Rh^{-1} = ({\cal R}_{\ty 21}^{-1})_{(V \ot V) \ot \Ch} =
                            a' \ty P \Rh\ty' + b' \ty P + c' \ty P K \;, \]
with some coefficients $a',b',c' \in \Ch\ty$. Since the tensor $B$ obviously
is symmetric, we conclude from Eqn.~\reff{Runiv} that
 \[ {\cal R}_{\ty 21}^{-1} = (1 \ot 1 + \ldots\,) \ty\ty e^{-h B} \;, \]
where the dots stand for terms similar to those in Eqn.~\reff{Runiv}.
Proceeding as above, we can show that
  \[ a' = 1 \quad,\quad c' = -b' = q - q^{-1} \;, \]
which implies that
  \[ \begin{array}{rcl}
  \Rh^{-1} \eqq \Rh\ty' - (q - q^{-1}) \ty \II \ot \II + (q - q^{-1}) \ty K
                                                                   \\[1.0ex]
           \eqq (q + q^{-1}) P_s -q \tty \II \ot \II
                                    + (q - q^{-1})(1 + q^{-2d-2})^{-1} K \;.
     \end{array} \]
It is easy to show directly that the operator on the right hand side really
is the inverse of $\Rh = \Rh\ty'$, which is a first check that our
calculations are correct. \vspace{2.0ex}

Summarizing part of the results of the present section, we have shown that
  \[ \begin{array}{rcl}
    \Rh &\! = \!& \dis{\sum_i \si_i \ty\tty q^{\si_i} \ty
                                               E_{i,\ty i} \ot E_{i,\ty i}
              \,+\, \sum_i \si_i \ty\tty q^{-\si_i}
                                  \ty E_{-i,\ty i} \ot E_{i,-i}}   \\[2.7ex]
        &   &\hspace{-1.2em}
             {} + \dis{\sum_{i \neq j,-j} \! \si_{i,\ty j} \ty\ty
                                      E_{j,\ty i} \ot E_{i,\ty j}} \\[2.7ex]
        &   &\hspace{-1.2em}
             {} + (q - q^{-1})\ty \dis{\sum_{i < j}
                                      E_{i,\ty i} \ot E_{j,\ty j}} \\[2.7ex]
        &   &\hspace{-1.2em}
             {} - (q - q^{-1})\ty \dis{\sum_{i < j}
              \Cqi_{-j,\ty j} \, C^q_{i,-i} \, E_{-j,\ty i} \ot E_{j,-i}} \;,
     \end{array} \]
or equivalently, that
  \[ \begin{array}{rcl}
      R &\! = \!& \dis{\sum_i q^{\si_i} \ty E_{i,\ty i} \ot E_{i,\ty i}
              \,+\, \sum_i q^{-\si_i}
                                  \ty E_{i,\ty i} \ot E_{-i,-i}}   \\[2.7ex]
        &   &\hspace{-1.2em}
             {} + \dis{\sum_{i \neq j,-j} \!
                                      E_{i,\ty i} \ot E_{j,\ty j}} \\[2.7ex]
        &   &\hspace{-1.2em}
             {} + (q - q^{-1})\ty \dis{\sum_{i < j}
                 \si_{i,\ty j} \ty\ty E_{j,\ty i} \ot E_{i,\ty j}} \\[2.7ex]
        &   &\hspace{-1.2em}
             {} - (q - q^{-1})\ty \dis{\sum_{i < j} \si_j \tty
              \Cqi_{-j,\ty j} \, C^q_{i,-i} \, E_{j,\ty i} \ot E_{-j,-i}} \;.
     \end{array} \]
In these equations, the indices $i$ and $j$ run through the index set $I$,
subject to conditions as specified. The equations hold in both settings,
the one in terms of formal power series, and the one where $q$ is a complex
number. It should also be noted that $E_{i,\ty j} \ot E_{k \ty,\ty\ty \ell}$
denotes the {\em normal} (non--super) tensor product of linear mappings (the
graded tensor product of two linear maps $f$ and $g$ is denoted by
$f \otb g\ty$).

\sect{Properties of $R$ and $\Rh$ \vspace{-1ex}}
In the present section, we want to collect some of the basic relations
satisfied by $R$ and $\Rh\ty\ty$. First of all, we recall the following
equations, which have been derived in the preceding section:
 \be \Rh \,=\, (q + q^{-1})\ty P_s - q^{-1} \II \ot \II
                        - (q - q^{-1})(1 + q^{2d+2})^{-1} K \label{RPs} \ee
 \be \Rh^{-1} \,=\, (q + q^{-1}) \ty P_s -q \ty\ty \II \ot \II
                    + (q - q^{-1})(1 + q^{-2d-2})^{-1} K \;.
                                            \hspace{0.9em} \label{RiPs} \ee
Using the results of Section 5 (in particular, Eqn.~\reff{bofa}), the first
of these equations implies that
  \[ \begin{array}{rcll}
      \Rh(x) \eqq q \tty x & \mbox{for all $x \in (V \ot V)_s$} \\[1.0ex]
      \Rh(y) \eqq -q^{-1} \ty y &
                           \mbox{for all $y \in (V \ot V)_a^0$} \\[1.0ex]
      \Rh(a) \eqq -q^{-2d-1} \ty a &.
     \end{array}\]
Moreover, the linear map induced by $\Rh$ in the one--dimensional
$\Uq$--module $(V \ot V)_a / (V \ot V)_a^0$ is equal to the multiplication by
$-q^{-2d-1}$. For $n \neq m\tty$, this follows from the last equation.

Next we recall that the $\Uq$--module endomorphisms of $V \ot V$ commute
one with another. In particular, we have
  \[ P_s \ty K = K P_s = 0 \]
  \[ \Rh\ty\ty K = K \Rh = -q^{-2d-1} K \;. \]
On the other hand, the Eqns.~\reff{RPs} and \reff{RiPs} imply that
  \[ \Rh -\Rh^{-1} = (q - q^{-1})(\II \ot \II - K) \;. \]
Obviously, this equation is equivalent with
  \[ \Rh^2 -(q - q^{-1}) \Rh - \II \ot \II
                                        = - (q - q^{-1}) \ty \Rh\ty K \;, \]
and hence also with
  \be (\Rh - q \ty \II \ot \II)(\Rh + q^{-1} \II \ot \II)
                   \,=\, (q - q^{-1}) \tty q^{-2d-1} K \;. \label{KofRh} \ee
Since the image of the operator $K$ is contained in $\CC\,a$\ty\ty, it
follows that
  \be (\Rh - q \ty \II \ot \II)(\Rh + q^{-1} \II \ot \II)
            (\Rh + q^{-2d-1} \ty \II \ot \II) \,=\, 0 \;. \label{minpol} \ee
It is easy to see that the polynomial $(X - q)(X + q^{-1})(X + q^{-2d-1})$
\vspace{0.2ex}
involved in Eqn.~\reff{minpol} is the minimal polynomial of the operator
$\Rh\ty\ty$.

The preceding equations can be used to write the spectral projectors of $\Rh$
as polynomials in $\Rh$ (to the extent in which these projectors exist). For
example, we find \vspace{-0.4ex}
 \[ P_s
    \,=\, \frac{(\Rh + q^{-1})(\Rh + q^{-2d-1})}{(q + q^{-1})(q + q^{-2d-1})}
                                                      \vspace{0.4ex} \;\,. \]
Moreover, we stress that according to Eqn.~\reff{KofRh}, the operator $K$ can
be written as a polynomial in $\Rh\ty\ty$. This fact (which is not true in
the undeformed case $q = 1$) will turn out to be crucial in the construction
of the quantum supergroup $\SP$.

Let us now derive two relations which are related to the fact that on $V$
there exists an invariant bilinear form, namely, the form $b^q$ found in
Section 4. We use the results of Appendix A and argue as in Ref.~\cite{KSc},
for the original setting in which $q$ is a complex number. The reader who is
not satisfied by this sloppy procedure may either reformulate everything
in terms of formal power series, or else regard the final result
Eqn.~\reff{Rfrl} as a conjecture which has to be (and has been) checked
independently.

Let
 \[ f_{\ell} : V \lra V^{\ast \rm gr} 
                                \quad,\quad f_r : V \lra V^{\ast \rm gr} \]
be the linear maps associated to $b^q$ (see Appendix A). Like $b^q$ they are
homogeneous of degree zero.

We write the universal $R$--matrix of $\Uq$ in the form
 \[ {\cal R} = \sum_\su R_s^1 \ot R_s^2 \;, \vso \]
where $R_s^1\ty,R_s^2 \in \Uq$. It is well--known that
 \[  {\cal R}^{-1} = (S \ot \id)({\cal R}) \;. \] 
This implies that
 \[ R^{-1} \,=\, {\cal R}_{V \ot V}^{-1}
           \,=\, (S \ot \id)({\cal R})_{V \ot V}
           \,=\, \sum_\su S(R_s^1)_V \otb (R_s^2)_V \;, \vso \]
where $\otb$ denotes the tensor product of linear maps {\em in the graded
sense.} Using Eqn.~\reff{linfr}, we conclude that
  \[ \begin{array}{rcl}
     R^{-1} \eq \dis{\sum_\su
                   f_r^{-1} \comp ((R_s^1)_V)^{\rm st} \comp f_r
                                                 \otb (R_s^2)_V } \\[-0.5ex]
            \eq \dis{(f_r^{-1} \ty\otb\ty\ty \II)
                      \comp \Big(\sum_\su \ty
                               ((R_s^1)_V)^{\rm st} \otb (R_s^2)_V \Big)
                                           \comp (f_r \otb \II) } \\[-0.5ex]
            \eq (f_r^{-1} \ot \II) \comp R^{\tty{\rm st}_1}
                                                \comp (f_r \ot \II) \;,
     \end{array} \]
where ${}^{{\rm st}_1}$ denotes the super--transposition of the first
tensorial factor (for more details, see Appendix B).

Similarly, we can start from the equation
   \[ {\cal R}^{-1} = (\id \ot S^{-1})({\cal R}) \]                         
and derive that
  \[ R^{-1} \,=\, {\cal R}_{V \ot V}^{-1}
                  \,=\, \sum_\su (R_s^1)_V \otb S^{-1}(R_s^2)_V \;. \vso \]
According to Eqn.~\reff{linfl}, this implies that
  \[ \begin{array}{rcl}
     R^{-1} \eq \dis{\sum_\su
                   (R_s^1)_V \otb
                   f_{\ell}^{-1} \comp ((R_s^2)_V)^{\rm st} \comp f_{\ell} }
                                                                  \\[-0.5ex] 
            \eq \dis{(\II \otb f_{\ell}^{-1})
                          \comp \Big(\sum_\su \ty
                            (R_s^1)_V \otb ((R_s^2)_V)^{\rm st}\ty \Big)
                                      \comp (\II \otb f_{\ell}) } \\[-0.5ex]
            \eq (\II \ot f_{\ell}^{-1}) \comp R^{\tty{\rm st}_2}
                                   \comp (\II \ot f_{\ell}) \;,
     \end{array} \]
where ${}^{{\rm st}_2}$ denotes the super--transposition of the second
tensorial factor. Summarizing, we have shown that
 \be R^{-1} \,=\, (f_r^{-1} \ot \II) \comp R^{\tty{\rm st}_1}
                                                   \comp (f_r \ot \II)
            \,=\, (\II \ot f_{\ell}^{-1}) \comp R^{\tty{\rm st}_2}
                              \comp (\II \ot f_{\ell}) \;. \label{Rfrl} \ee
Note that these equations imply that $R^{\tty{\rm st}_1}$ and
$R^{\tty{\rm st}_2}$ are invertible.

The equations \reff{Rfrl} can be checked directly. To do that we need the
matrices of the linear maps $f_{\ell}$ and $f_r\ty\ty$. If $(e'_i)_{i \in I}$
is the basis of $V^{\ast\rm gr}$ dual to $(e_i)_{i \in I}\ty\ty$, we find for
all $j \in I$
  \[ f_{\ell}(e_j) = \sum_{i \in I} \ty C_{j,\ty i}^q \, e'_i
                                     \quad,\quad
          f_r(e_j) = \sum_{i \in I} \si_{j,\ty i}
                                         \ty\ty C_{i,\ty j}^q \, e'_i \;. \]
We also need a formula for $R^{-1}$. Using the equation
  \[ R^{-1} \,=\, \Rh^{-1} P \,=\, PR\ty\tty P - (q - q^{-1})(P - KP) \;, \]
we derive that
  \[ \begin{array}{rcl}
     R^{-1} &\! = \!& \dis{\sum_i q^{-\si_i} \ty E_{i,\ty i} \ot E_{i,\ty i}
              \,+\, \sum_i q^{\si_i}
                                  \ty E_{i,\ty i} \ot E_{-i,-i}}   \\[2.7ex]
        &   &\hspace{-1.2em}
             {} + \dis{\sum_{i \neq j,-j} \!
                                      E_{i,\ty i} \ot E_{j,\ty j}} \\[2.7ex]
        &   &\hspace{-1.2em}
             {} - (q - q^{-1})\ty \dis{\sum_{i < j}
                 \si_{i,\ty j} \ty\ty E_{j,\ty i} \ot E_{i,\ty j}} \\[2.7ex]
        &   &\hspace{-1.2em}
             {} + (q - q^{-1})\ty \dis{\sum_{i < j} \si_i \tty
              \Cqi_{j,-j} \, C^q_{-i,\ty i} \, E_{j,\ty i} \ot E_{-j,-i}} \;.
     \end{array} \]
Recalling the formulae for the partial super--transpose given in Appendix B,
it is now not difficult to show that the equations \reff{Rfrl} are indeed
satisfied.

A closer look at the formula for $R^{-1}$ reveals that, somewhat
unexpectedly, $R_q^{-1}$ is {\em not} equal to $R_{q^{-1}}$ (we are using
the obvious notation). This fact is closely related to Eqn.~\reff{Cqsq}.

In the purely symplectic case considered in Ref.~\cite{RTF} it is known that
$\Rh^{\tty\rm t_1 t_2}$ is equal to $\Rh$ (where ${}^{\ty\rm t_1}$ and
${}^{\ty\rm t_2}$ denote the usual transposition of the first resp.~second
tensorial factor). For reasons similar to those above, I have not been able
to derive an analogous equation in the present setting.

We proceed by recalling that the general theory of quasitriangular Hopf
superalgebras implies that $R$ satisfies the {\em graded}\/ Yang--Baxter
equation. Equivalently, this means that $\Rh$ satisfies the braid relation
 \[ (\Rh \ot \II)(\II \ot \Rh)(\Rh \ot \II)
                       \,=\, (\II \ot \Rh)(\Rh \ot \II)(\II \ot \Rh) \;. \]
It would be worth--while to check this relation directly, but I haven't done
that.

Finally, I have shown by explicit calculation that $\Rh$ and $K$ satisfy the
following relations:
 \[ (\II \ot K)(\Rh^{\pm 1} \ot \II)(\II \ot K)
                                      \,=\, - q^{\pm(2d+1)}(\II \ot K) \]
 \[ \;\;
    (K \ot \II)(\II \ot \Rh^{\pm 1})(K \ot \II)
                                      \,=\, - q^{\pm(2d+1)}(K \ot \II) \;. \]
Summarizing part of the results of the present section, we conclude that
$\Rh$ and $K$ generate representations of the Birman, Wenzl, Murakami
algebras \cite{BWe}, \cite{Mur} as defined in Ref.~\cite{LRa} (with
$z = -q^{2d+1}$).

\sect{Comparison with known special cases \vspace{-1ex}}
In a few special cases, the $R$--matrix calculated in this work has already
been known.

\noindent
1. The case $m = 0$ \\[1.0ex]
It should be obvious to the reader that our results apply in the case
$m = 0$ as well, and this case has been settled in Ref.~\cite{RTF}. Actually,
I have used this fact throughout the whole investigation: It enabled me to
check my calculations and to guess a concise expression for the $R$--matrix.
For greater clarity, let us mark the entries of the present work by the
subscript ``here'' and those of Ref.~\cite{RTF} by the subscript ``RTF'',
moreover, let us indicate the dependence on the parameter $q$ by the
superscript $q\ty\ty$. Then we have
 \[ R^q_{\rm here} = R^q_{\rm RTF} \quad,\quad 
                                  \Rh^q_{\rm here} = \Rh^q_{\rm RTF} \;, \]
and also
  \[ C^q_{\rm here} = -C^q_{\rm RTF} \quad,\quad
                                     K^q_{\rm here} = K^q_{\rm RTF}  \;. \]

\noindent
2. The case $n = 0$ \\[1.0ex]
This case is more interesting. Once again, this case has been treated in
Ref.~\cite{RTF}. On the other hand, the calculations of the present work
don't make sense in this case from the outset, since the root system of the
Lie algebra ${\mathfrak o}(2m)$ does not have a basis of the type used here.
Nevertheless, we find
 \[ R^{\ty q}_{\rm here} = R^{(q^{-1})}_{\rm RTF} \quad,\quad 
                  \Rh^{\ty q}_{\rm here} = -\Rh^{(q^{-1})}_{\rm RTF} \;, \]
furthermore, we have
  \[ C^q_{\rm here} = q^{-1} \ty C^{(q^{-1})}_{\rm RTF} \quad,\quad
                            K^q_{\rm here} = K^{(q^{-1})}_{\rm RTF}  \;. \]
The change from $q$ to $q^{-1}$ under the transition from the even to the odd
case is a known phenomenon. On the other hand, the sign factor in the formula
for $\Rh$ is easily understood: In the purely odd case, the supersymmetric
twist is equal to minus the normal (non--graded) twist. \vspace{1.0ex}

\noindent
3. The case $m = n = 1$ \\[1.0ex]
In this case, the universal $R$--matrix and the $R$--matrix in the vector
representation have been given in Ref.~\cite{DFI}. However, these authers
have worked with a basis of the root system which consists of two odd roots
(see Ref.~\cite{ACF}). Actually, I have applied the approach of the present
paper also to this case, and (after some obvious adjustments) indeed have
obtained the $R$--matrix of Ref.~\cite{DFI}.

\sect{Discussion \vspace{-1ex}}
We have calculated the $R$--matrix of the symplecto--orthogonal quantum
superalgebra $\Uq$ in the vector representation, and we have derived its most
important properties. In a subsequent work \cite{Sqg}, we shall use this
$R$--matrix to construct the corresponding quantum supergroup $\SP$ and its
basic comodule superalgebras.

A special feature of the present work is that we have used a somewhat
unusual basis of the root system of $\spo$. This was dictated by the wish for
a unified treatment of all cases, and by the assumption that the basis of
the root system should contain only one odd root. If one drops this last
requirement, there is another possibility: For $m \geq 2\tty$, one chooses
Kac's distinguished basis, whereas for $m = 1$ (i.e., for the $C$--type Lie
superalgebras) one chooses a basis containing two odd simple roots. In the
latter case, the Dynkin diagram looks as follows:
\begin{center}
\unitlength1.25mm
\vspace*{-11ex}
\begin{picture}(80.2,30)
\thicklines
\put(2.2,5){\circle{4.4}}
\put(4.4,5){\line(1,0){10}}
\put(16.6,5){\circle{4.4}}
\put(18.8,5){\line(1,0){10}}
\multiput(30.3,4.8)(2.0,0){4}{.}
\put(38.8,5){\line(1,0){10}}
\put(51.0,5){\circle{4.4}}
\put(53.2,5){\line(1,0){10}}
\put(65.4,5){\circle{4.4}}
\put(67.5,6.0){\line(2,1){9.8}}
\put(79.4,11.9){\circle{4.4}}
\put(78.0,10.4){\line(1,1){2.9}}
\put(78.0,13.4){\line(1,-1){2.9}}
\put(67.5,4.0){\line(2,-1){9.8}}
\put(79.4,-1.9){\circle{4.4}}
\put(78.0,-3.4){\line(1,1){2.9}}
\put(78.0,-0.4){\line(1,-1){2.9}}
\put(78.9,0.3){\line(0,1){9.4}}
\put(79.9,0.3){\line(0,1){9.4}}
\end{picture}
\end{center}
\vspace{2.1ex}
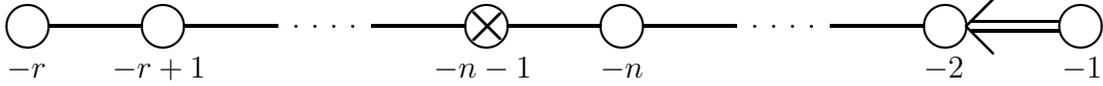
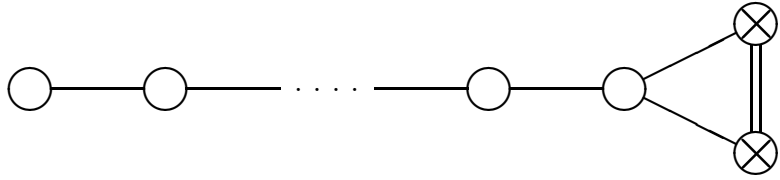
\centerline{Dynkin diagram of the Lie superalgebra
                                            $\mathfrak{spo}(2n \ty|\ty 2)$}
\vspace{2.5ex}

\noindent
For $n = m = 1\ty$, this is just the choice mentioned in the preceding
section. It should be interesting to calculate the $R$--matrix also under
these assumptions.

In Section 6 we had to use a formulation of the theory in terms of formal
power series. This made our arguments somewhat clumsy. Of course, we could
exclusively use the language of formal power series. The present formulation
has been chosen in view of possible applications.
\vspace{4.0ex}

\noindent
{\LARGE\bf Appendix \vspace{-1.5ex}}

\begin{appendix}

\sect{Invariant bilinear forms \vspace{-1ex}}
In the following, the base field will be an arbitrary field $\KK$ of
characteristic zero. Let $\Ga$ be an abelian group, and let $\si$ be a
commutation factor on $\Ga$ with values in $\KK\,$. It is well-known that
the class of $\Ga$--graded vector spaces, endowed with the usual tensor
product of graded vector spaces and with the twist maps defined by means of
$\si\tty$, forms a tensor category (see Ref.~\cite{Sgt} for details). A
(generalized) Hopf algebra $H$ living in this category is called a
$\si$--Hopf algebra. More explicitly, $H$ is a $\Ga$--graded associative
algebra with a unit element, and it is endowed with a coproduct $\De\ty\ty$,
a counit $\ve\ty\ty$, and an antipode $S\ty\ty$, which satisfy the obvious
axioms (in the category). In particular, this implies that the structure
maps $\De\ty\ty$, $\ve\ty\ty$, and $S$ are homogeneous of degree zero. In
the following, we shall freely use the notation and results of
Ref.~\cite{Sgt}. (For generalized Hopf algebras living in more general
categories, see Ref.~\cite{Maj}.)

Let $V$ and $W$ be two graded (left) $H$--modules. Then $V \ot W$ and
$\Lgr(V,W)$ have a canonical structure of a graded $H$--module as well.
(Recall that $\Lgr(V,W)$ denotes the space of all linear maps of $V$ into $W$
which can be written as a sum of homogeneous linear maps of $V$ into $W$.)
If $U$ is a third graded $H$--module, there exists a canonical isomorphism
of graded $H$--modules,
 \be \lam : \Lgr(V,\Lgr(W,U)) \lra \Lgr(V \ot W,U) \;, \label{lam} \ee
which is defined by
 \[ \lam(f)(x \ot y) = (f(x))(y) \;, \]
for all $f \in \Lgr(V,\Lgr(W,U))$, $x \in V$, and $y \in W$.

The next thing to be mentioned is that an element $x$ of a graded $H$--module
$V$ is said to be $H$--invariant (or simply invariant) if
 \[ h \cd x = \ve(h) \ty\ty x \quad\mbox{for all $h \in H$} \,. \]
(Recall that, quite generally, the dot denotes a module action.) Let
$g \in \Lgr(V,W)$ be homogeneous of degree $\ga$. Then $g$ is invariant if
and only if it is $H$--linear in the graded sense, i.e., if and only if
 \[ g(h \cd x) = \si(\ga,\eta) \ty\ty h \cd g(x) \;, \]
for all elements $h \in H_{\eta}\ty\ty$, $\eta \in \Ga$, and all $x \in V$.
 
In the following, we choose $U = \KK\,$, the trivial $H$--module. Then 
$\Lgr(V \ot W,\KK\ty) = (V \ot W)^{\ast\rm gr}$ is the graded dual of
$V \ot W$. It is well--known that, regarded as a graded vector space, this
space is canonically isomorphic to $\Lgr_2(V,W\ty;\KK\ty)$, the space of all
bilinear forms on $V \ti W$ that can be written as a sum of homogeneous
bilinear forms on $V \ti W$. The canonical isomorphism is used to transfer
the $H$--module structure from $\Lgr(V \ot W,\KK\ty)$ to
$\Lgr_2(V,W\ty;\KK\ty)$. For every bilinear form
$b \in \Lgr_2(V,W\ty;\KK\ty)$, the corresponding linear form on $V \ot W$
will be denoted by $\bt\ty\ty$.

Now let
  \[ b : V \ti W \lra \KK \]
be a bilinear form on $V \ti W$ that is homogeneous of degree $\bet$\ty, let
  \[ \bt : V \ot W \lra \KK \]
be the associated linear form on $V \ot W$, and let
  \[ f_{\ell} : V \lra W^{\ast \rm gr} \]
be the linear map canonically corresponding to $\bt\ty\ty$, i.e.,
  \[ (f_{\ell}(x))(y) = b(x,y) \]
for all $x \in V$ and $y \in W$. (Choosing $U = \KK$ in Eqn.~\reff{lam},
this is to say that $\lam(f_{\ell}) = \bt\ty\ty$.) Then $\bt$ and $f_{\ell}$
are homogeneous of degree $\bet$\ty. According to the preceding discussion,
the following statements are equivalent: \\[1.0ex]
1) The bilinear form $b$\,, or equivalently, the linear form $\bt$\,, is
$H$--invariant, i.e., we have
  \[ \si(\eta,\bet)\, \bt \comp S(h)_{V \ot W} = \ve(h)\ty \bt \] 
for all elements $h \in H_{\eta}\ty\ty$, $\eta \in \Ga$.
                                                                     \\[1ex]
2) The linear form $\bt$ is $H$--linear in the graded sense, i.e., we have
  \[ \bt\ty(h \cd (x \ot y))
                           = \si(\bet,\eta)\ty \ve(h)\ty \bt\ty(x \ot y) \]
for all elements $h \in H_{\eta}\ty\ty$, $\eta \in \Ga$, and all $x \in V$,
$y \in W$.  Since $\ve(h) = 0$ if $\eta \neq 0\tty$, this
is equivalent with
 \be \bt\ty(h \cd (x \ot y)) = \ve(h)\ty \bt\ty(x \ot y) \label{linbt} \ee
for all elements $h \in H$, $x \in V$, and  $y \in W$. \\[1ex]
3) The linear map $f_{\ell}$ is $H$--invariant. \\[1ex]
4) The linear map $f_{\ell}$ is $H$--linear in the graded sense, i.e., we
have
  \[ f_{\ell}\ty(h \cd x) = \si(\bet,\eta)\tty h \cd f_{\ell}(x) \]
for all elements $h \in H_{\eta}\ty\ty$, $\eta \in \Ga$, and all $x \in V$.
This is equivalent with 
  \be f_{\ell} \comp h_V = \si(\bet,\eta)(S(h)_W)^{\rm T} \comp f_{\ell}
                                                          \label{linfl} \ee
for all elements $h \in H_{\eta}\ty\ty$, $\eta \in \Ga$, and also with
 \be b(h \cd x\ty,y) = \si(\eta\ty,\xi)\tty b(x\ty,S(h) \cd y) \label{binv} \ee
for all elements $h \in H_{\eta}\ty\ty$, $\eta \in \Ga$, all
$x \in V_{\xi}\ty\ty$, $\xi \in \Ga$, and all $y \in W$. Recall that
${}^{\rm T}$ denotes the $\si$--transposition. (In the super case, the
super--transposition will be denoted by ${}^{\rm st}$.)

To proceed, we note that, apart from $f_{\ell}\ty\ty$, the bilinear form $b$
yields a second linear map, namely, the map 
  \[ f_r : W \lra V^{\ast \rm gr}  \]
which, for all elements $x \in V_{\xi}\ty\ty$, $\xi \in \Ga$, and
$y \in W_{\eta}\ty\ty$, $\eta \in \Ga$, is given by
  \[ (f_r(y))(x) = \si(\eta\ty,\xi)\ty b(x,y) \;. \]
If
  \[ \nu_V : V \lra (V^{\ast \rm gr})^{\ast \rm gr} \quad,\quad
                         \nu_W : W \lra (W^{\ast \rm gr})^{\ast \rm gr} \]
are the canonical injections (in the graded sense, see Ref.~\cite{Sgt}), we
have
  \[ f_r = f_{\ell}^{\rm T} \comp \nu_W \quad,\quad
                                 f_{\ell} = f_r^{\rm T} \comp \nu_V \;. \]
Then the condition \reff{linfl} is equivalent with
  \be f_r \comp S(h)_W = \si(\bet,\eta) (h_V)^{\rm T} \comp f_r
                                                      \;, \label{linfr} \ee
for all elements $h \in H_{\eta}\ty\ty$, $\eta \in \Ga$.
Indeed, the Eqns.~\reff{linfl} and \reff{linfr} can be derived from each
other by composing their $\si$--transposes with $\nu_W$ and $\nu_V\ty\ty$,
respectively.

\sect{Partial transposition \vspace{-1ex}}
In the following, the base field will be an arbitrary field $\KK$ of
characteristic zero. Let $\Ga$ be an abelian group, and let $\si$ be a
commutation factor on $\Ga$ with values in $\KK\,$. All gradations considered
in this appendix will be $\Ga$--gradations. We shall freely use the notation
and results of Ref.~\cite{Sgt}.

Let $V$ and $W$ be two {\em finite--dimensional} graded vector spaces. Then
there exists a unique linear map
  \[ {}^{{\rm T}_1} : \Lgr(V \ot W,V \ot W)
              \lra \Lgr(V^{\ast \rm gr} \ot W,V^{\ast \rm gr} \ot W) \;, \]
such that
  \[ (f \,\otb\, g)^{{\rm T}_1} = f^{\rm T} \,\otb\, g  \;, \]
and a unique linear map
  \[ {}^{{\rm T}_2} : \Lgr(V \ot W,V \ot W)
              \lra \Lgr(V \ot W^{\ast \rm gr},V \ot W^{\ast \rm gr}) \;, \]
such that
  \[ (f \,\otb\, g)^{{\rm T}_2} = f \,\otb\, g^{\rm T} \;,\]
for all $f \in \Lgr(V,V)$ and all $g \in \Lgr(W,W)$. (Recall that $\otb$
denotes the graded tensor product of linear maps, and that ${}^{\rm T}$
denotes the $\si$--transposition.) We call ${}^{{\rm T}_1}$ and
${}^{{\rm T}_2}$ the partial $\si$--transposition of the first and second
tensorial factor, respectively. (In the super case, we shall write
${}^{{\rm st}_1}$ and ${}^{{\rm st}_2}$.)

It is easy to see that
  \[ (f \ot g)^{{\rm T}_2} = f \ot g^{\rm T} \;, \]
for all $f \in \Lgr(V,V)$ and all $g \in \Lgr(W,W)$. On the other hand, under
our present general assumptions, $(f \ot g)^{{\rm T}_1}$ {\em is not,} in
general, equal to $f^{\rm T} \ot g \ty\ty$. Using our standard notation for
$V$, and a similar notation for $W$, but with all entries overlined, we find
instead that
  \[ (E_{ij} \ot \ol{E}_{rs})^{{\rm T}_1}
           = \si(\eta_i + \eta_j,\ol{\eta}_r - \ol{\eta}_s) \ty
                                     E_{ij}^{\rm T} \ot \ol{E}_{rs} \;. \] 
{\em In the super case,} which is the case we are mainly interested in, this
equation implies that
  \[ (f \ot g)^{{\rm T}_1} = \si(\vp,\ga) f^{\rm T} \ot g \;, \]
where $f \in \Lgr(V,V)$ is homogeneous of degree $\vp$ and $g \in \Lgr(W,W)$
is homogeneous of degree $\ga$.

\end{appendix}

\end{document}